\documentclass[10pt]{article}

\usepackage{a4wide}
\usepackage{amssymb}
\usepackage{amsfonts}
\usepackage{amsmath}
\input xy
\xyoption{arrow} \xyoption{matrix}

\date{}

\newtheorem{proposition}{Proposition}[section]
\newtheorem{theorem}[proposition]{Theorem}
\newtheorem{lemma}[proposition]{Lemma}

\newtheorem{corollary}[proposition]{Corollary}

\def\GK{{\rm  GK}\,}

\def\der{\partial }

\def\nFM0{{\nu }_{F,M_0}}
\def\nFN0{{\nu }_{F,N_0}}
\def\nGN0{{\nu }_{G,N_0}}

\def\N0{ {\bf N}_0 }

\def\t{\otimes}

\def\v{\varphi}
\def\ra{\rightarrow}

\def\lra{\leftrightarrow}
\def\Xpm{X^{\pm }}

\def\s{\sigma}
\def\Z{\mathbb{Z}}

\def\l1{{\lambda}_1}

\def\a{\alpha}
\def\a0{ {\alpha }_0}
\def\a1{ {\alpha }_1}

\def\l{\lambda}

\def\nFGM0{{\nu }_{F,G,M_0}}

\def\nFN0{{\nu}_{F,N_0}}

\def\sm{{\sigma}^m}

\def\sm1{{\sigma}^{-1}}

\def\smtp1{{\sigma}^{-t+1}}

\def\S1{S^{-1}}

\def\Xpm1{X^{\pm 1}_1}

\def\sPM1{{\sigma }^{\pm 1}}
\def\sMP1{{\sigma }^{\mp 1 }}

\def\d{\delta}

\def\fd{{\rm fd}_A}

\def\di{{\rm d.ind}}

\def\L{\Lambda}

\def\OO{{\cal O}}
\def\CA{{\cal A}}

\def\Ytm1{Y^{t-1}}
\def\Yim1{Y^{i-1}}

\def\CK{{\cal K}}

\def\CS{{\cal S}}

\def\ass{{\rm ass}}

\def\supp{{\rm supp }}

\def\Aut{{\rm Aut}}

\def\dim{{\rm dim }}

\def\ker{ {\rm ker } }

\def\CJ{ {\cal J}}

\def\SL2Z{ {\rm SL}_2({\bf Z}) }

\def\th{ \theta }

\def\Gp1{ G^{1 , 1 } }
\def\P11{ P^{-1 , 1 } }
\def\Pp1{ P^{1 , 1 } }

\def\th{\theta}

\def\nCLsr{{}^\nu\kern-2pt {\cal L}^{\sigma , \rho  }}
\def\nP{{}^\nu \kern-2pt P}
\def\nL{{}^\nu\kern-2pt L}
\def\nLL{{}^\nu\kern-2pt \Lambda}
\def\nPsr{{}^\nu\kern-2pt P^{\sigma , \rho  }}
\def\nLsr{{}^\nu\kern-2pt L^{\sigma , \rho  }}
\def\nuCL{{}^\nu\kern-2pt  {\cal L}}
\def\nCLsr{{}^\nu\kern-2pt {\cal L}^{\sigma , \rho  }}
\def\nCL1m{{}^\nu\kern-2pt {\cal L}^{-1 , 1  }}
\def\x1nu{x^\frac{1}{\nu}}
\def\xm1nu{x^{-\frac{1}{\nu}}}

\def\bL{\bar{L}}

\def\ra{\rightarrow }

\def\CB{{\cal B}}

\def\CI{{\cal I}}

\def\CT{{\cal T}}

\def\CC{ {\cal C}}

\def\CP{ {\cal P}}

\def\nAM0{{\nu }_{{\cal A},M_0}}
\def\nAN0{{\nu }_{{\cal A},N_0}}

\def\End{ {\rm End }}

\def\CJ{ {\cal J }}

\def\CP{ {\cal P }}

\def\bI{\overline{I}}

\def\ga{\mathfrak{a}}
\def\gb{\mathfrak{b}}

\def\gd{\mathfrak{d}}

\def\gp{\mathfrak{p}}
\def\gq{\mathfrak{q}}

\def\SL{{\rm SL}}

\def\Spec{{\rm Spec}}

\def\di!{\frac{\der^i}{i!}}
\def\dik!{\frac{\der^k_i}{k!}}
\def\hA{\widehat{A}}

\def\Max{{\rm Max}}

\def\gldim{{\rm gldim}}

\def\N{\mathbb{N}}

\def\0{\overline{0}}
\def\1{\overline{1}}

\def\Ln1{\L_{n,\overline{1}}}

\def\a1{a_{\overline{1}}}

\def\S{\Sigma}

\def\vn1{\overrightarrow{n-1}}

\def\im{{\rm im}}

\def\mA{\mathbb{A}}

\def\soc{{\rm soc}}
\def\Sub{{\rm Sub}}

\def\Inn{{\rm Inn}}

\def\mS{\mathbb{S}}
\def\mJ{\mathbb{J}}
\def\mI{\mathbb{I}}

\def\clKdim{{\rm cl.Kdim}}
\def\pd{{\rm pd}}

\def\lgldim{{\rm l.gldim}}
\def\rgldim{{\rm r.gldim}}

\def\fd{{\rm fd}}

\def\mT{\mathbb{T}}

\def\mE{\mathbb{E}}

\def\Frac{{\rm Frac}}

\def\wdim{{\rm wdim}}

\def\K1{{\rm K}_1}

\def\wdim{{\rm w.dim}}

\def\hmI1{\widehat{\mI_1}}
\def\tmI1{\widetilde{\mI_1}}
\def\tmJ1{\widetilde{\mJ_1}}
\def\hB1{\widehat{B_1}}
\def\hCB1{\widehat{\CB_1}}

\def\ga{\mathfrak{a}}

\def\tor{{\rm tor}}

\def\pdIn{{\rm pd}_{\mI_n}}

\begin{document}

\author{V. V. \  Bavula 
}

\title{The global dimension of the  algebras of polynomial
integro-differential operators $\mI_n$ and the Jacobian algebras $\mA_n$}

\maketitle

\begin{abstract}
The aim of the paper is to prove two conjectures from the paper \cite{algintdif} that the (left and right)  global dimension of the algebra $\mI_n:=K\langle x_1, \ldots , x_n,
\frac{\der}{\der x_1}, \ldots ,\frac{\der}{\der x_n},  \int_1,
\ldots , \int_n\rangle $ of  polynomial   integro-differential operators  and the Jacobian algebra $\mA_n$ is equal to $n$ (over a field of characteristic zero). The algebras $\mI_n$ and $\mA_n$ are  neither left nor right Noetherian and $\mI_n\subset \mA_n$. Furthermore, they contain infinite direct sums of nonzero left/right ideals and are not domains. An analogue of  Hilbert's Syzygy Theorem is proven for the algebras $\mI_n$, $\mA_n$ and their factor algebras. It is proven that the  global dimension of all prime factor algebras  of the algebras  $\mI_n$  and  $\mA_n$ is $n$ and the weak global dimension of all the factor algebras of $\mI_n$ and $\mA_n$ is $n$.\\



$${\bf Contents}$$
\begin{enumerate}
\item Introduction.
\item The  global dimension of the algebra $\mI_1$.
 \item The  global dimension of the algebra $\mI_n$. \item The global dimension of the Jacobian algebra $\mA_n$.
 \item The weak  global dimension of  factor algebras of $\mI_n$.
\item The weak  global dimension of  factor algebras of
$\mA_n$.

\end{enumerate}
\end{abstract}


\section{Introduction}\label{INTRO}
Throughout,  $K$ is a
field of characteristic zero and  $K^*$ is its group of units; algebra  means an associative $K$-algebra  with $1$; module means
a left module;
 $\N :=\{0, 1, \ldots \}$ is the set of natural numbers;
$P_n:= K[x_1, \ldots , x_n]$ is a polynomial algebra over $K$;
$\der_1:=\frac{\der}{\der x_1}, \ldots , \der_n:=\frac{\der}{\der
x_n}$ are the partial derivatives ($K$-linear derivations) of
$P_n$; $\End_K(P_n)$ is the algebra of all $K$-linear maps from
$P_n$ to $P_n$;
the
subalgebra  $A_n:= K \langle x_1, \ldots , x_n , \der_1, \ldots ,
\der_n\rangle$ of $\End_K(P_n)$ is called the $n$'th {\em Weyl}
algebra.


{\it Definition}, \cite{Bav-Jacalg}. The {\em Jacobian algebra}
$\mA_n= K \langle x_1, \ldots , x_n , \der_1, \ldots ,
\der_n, H_1^{-1}, \ldots , H_n^{-1}\rangle$ is the subalgebra of $\End_K(P_n)$ generated by the Weyl
algebra $A_n$ and the elements $H_1^{-1}, \ldots , H_n^{-1}\in
\End_K(P_n)$ where $$H_1:= \der_1x_1, \ldots , H_n:= \der_nx_n.$$

The Jacobian algebras appeared in study of automorphisms of polynomial algebras as a very effective computational tool (which is not surprising  as they contain the algebras of polynomial integro-differential operators $\mI_n$, see below).  Clearly, $\mA_n =\bigotimes_{i=1}^n \mA_1(i) \simeq \mA_1^{\t n }$
where $\mA_1(i) := K\langle x_i, \der_i , H_i^{-1}\rangle \simeq
\mA_1$.  The algebra $\mA_n$ is neither a left nor right localization of the Weyl algebra $A_n$ at the multiplicative set generated by the elements $H_1, \ldots , H_n$ since the algebra $\mA_n$ is not a domain but the Weyl algebra $A_n$ is a domain.  The algebra $\mA_n$ contains all the  integrations
$\int_i: P_n\ra P_n$, $ p\mapsto \int p \, dx_i$ since  $$\int_i=
x_iH_i^{-1}: \; x_1^{\alpha_1} \cdots x_n^{\alpha_n}\mapsto (\alpha_i+1)^{-1}x_ix_1^{\alpha_1} \cdots x_n^{\alpha_n}\;\; {\rm where}\;\; \alpha_1,  \ldots ,\alpha_n\in \N .$$ In
particular, the algebra $\mA_n$ contains the {\em algebra}
$\mI_n:=K\langle x_1, \ldots , x_n$,
 $\der_1, \ldots ,\der_n,  \int_1,
\ldots , \int_n\rangle $ {\em of  polynomial integro-differential
operators}. Notice that $\mI_n=\bigotimes_{i=1}^n\mI_1(i)\simeq
\mI_1^{\t n}$ where $\mI_1(i):= K\langle x_i, \der_i,
\int_i\rangle$ and
$$P_n\subset A_n \subset \mI_n \subset \mA_n.$$

The algebras $\mI_n$ and $\mA_n$ are  neither left nor right Noetherian and  not domains.
 The algebras $\mI_n$ were introduced and  studied in detail in \cite{algintdif} and they have remarkable properties, name just a few (see  \cite{algintdif} for details):
\begin{itemize}
\item $\mI_n$ is a prime, catenary, central  algebra of classical Krull
dimension $n$ and of Gelfand-Kirillov
dimension $2n$, and there is a unique maximal ideal $\ga_n$ of the
algebra $\mI_n$.
\item  A canonical form is found for each element of
$\mI_n$ by showing that the algebra $\mI_n$ is a generalized Weyl
algebra.
\item The lattice $\CJ (\mI_n)$  of ideals of the algebra $\mI_n$ is a distributive lattice.
\item $\ga \gb=\gb\ga $ and $\ga^2= \ga$ for all
$\ga , \gb \in \CJ (\mI_n)$.   \item Classifications of all the ideals and the
prime ideals of the algebra $\mI_n$ are given. \item The set of ideals $\CJ
(\mI_n)$ of $\mI_n$ is a finite set and the number  $s_n:= |\CJ (\mI_n)|$ is equal to the {\em Dedekind number} ($2-n+\sum_{i=1}^n2^{n\choose i}\leq s_n
\leq 2^{2^n}$). \item $P_n$ is the only
(up to isomorphism) faithful simple $\mI_n$-module.
\item  Two sets of
defining relations are given for the algebra $\mI_n$.
\item The factor algebra
$\mI_n/ \ga$ is a  Noetherian algebra iff the ideal $\ga$ is equal to the unique  maximal ideal of $\mI_n$.
\item $\GK (\mI_n/\ga)=2n$ for all
ideals $\ga$ of $\mI_n$ distinct from $\mI_n$ (where $\GK$ is the Gelfand-Kirillov dimension).
\item The algebra $\mI_n$ admits an involution $*$ given by the rule
 $\der_i^*=\int_i$, $\int_i^* = \der_i$,
and $H_i^* = H_i$.
\item $\ga^* = \ga$ for all ideals $\ga$ of the algebra
$\mI_n$.
\item Each ideal of the algebra
$\mI_n$ is an essential left and right submodule of $\mI_n$.
\end{itemize}

The fact  that   certain rings of differential operators are catenary
was proven   by  Brown, Goodearl and  Lenagan in
\cite{Brown-Goodearl-Lenagan}.

The group of automorphisms of the algebra $\mI_n$ was found in \cite{Intdifaut}, it is an iterated semi-direct product of three obvious subgroups \cite[Theorem 5.5.(1)]{Intdifaut}.
\begin{itemize}
\item (\cite[Theorem 5.5.(1)]{Intdifaut}) $\; \Aut_K(\mI_n)=S_n\ltimes \mT^n\ltimes \Inn (\mI_n)$ {\em where $S_n$ is the symmetric group (it permutes the tensor components of the algebra $\mI_n = \mI_1^{\t n}$), $\mT^n$ is the $n$-dimensional algebraic torus and $\Inn (\mI_n)$ is the group of inner automorphisms of the algebra} $\mI_n$.
\end{itemize}

In \cite{Intdifaut}, the group of units $\mI_n^\times$ of the algebra $\mI_n$ is described.
\begin{itemize}
\item  (\cite[Theorem 3.1.(2)]{IntdifDixConj}) $\mI_n^\times = K^*\times (1+\ga_n)^\times$ {\em where $\ga_n$ is the unique maximal ideal of the algebra $\mI_n$ and $(1+\ga_n)^\times$ is the group of units of the multiplicative monoid $1+\ga_n$. }
\end{itemize}

The group $\Inn (\mI_n)\simeq (1+\ga_n)^\times$ is huge.
 The groups of automorphisms of the polynomial algebras $P_n$ (resp., the Weyl algebras $A_n$) are found only for $n=1,2$ (resp., $n=1$, see \cite{Dix}). The groups of automorphisms of the algebras $P_2$, $A_1$ and $\mI_1$ have similar structure. One of the serious obstacles in finding the groups of automorphisms for the polynomial and Weyl algebras are the {\em Jacobian Conjecture} (for polynomials) and the {\em Dixmier Conjecture} (for the Weyl algebras) that are still open. The Dixmier Conjecture states that {\em every algebra endomorphism of the Weyl algebra $A_n$ is an automorphism, \cite{Dix}.} In \cite{IntdifDixConj}, it is shown that an analogue of the Dixmier Conjecture holds for the algebra $\mI_1$:

\begin{itemize}
\item  (\cite[Theorem 1.1]{IntdifDixConj}) {\em Every algebra endomorphism of the algebra $\mI_1$ is an automorphism.}
\end{itemize}

In \cite{IntdifDixConj}, it was conjectured that the same is true for all algebras $\mI_n$. The algebras $\mI_n$ and $\mA_n$ are closely related and have similar properties, see \cite{Bav-Jacalg,AutJacAlg}.

{\bf The global dimension of the algebras $\mI_n$, $\mA_n$ and their prime factor algebras.}  In 1962, Rinehart found the global dimension of the first Weyl algebra $A_1$, \cite{Reinhart-gldim-A1}.  In 1972, Roos proved that the global dimension of
the Weyl algebra $A_n$ is $n$, \cite{Roos}. Goodearl obtained
formulae for the global dimension of certain rings of differential
operators  \cite{Goodearl-gldim-II}, \cite{Goodearl-gldim-III}. The weak global  dimension of the algebras $\mI_n$ and $\mA_n$ is $n$, \cite[Theorem 6.2, Theorem 7.2]{algintdif}. Since the algebras $\mI_n$ and $\mA_n$ admit an involution, they are {\em self-dual algebras} (i.e., they are isomorphic to their {\em opposite algebras}). As a result, for the algebras $\mI_n$ and $\mA_n$,  the left global dimension ($\lgldim$) coincides with the right global dimension ($\rgldim$) and we use the notation $\gldim$ for their common value.
\begin{itemize}
\item  (\cite[Proposition 6.7]{algintdif}) $n\leq \gldim (\mI_n)\leq 2n$.
\item  (\cite[Proposition 7.5]{algintdif}) $n\leq \gldim (\mA_n)\leq 2n$.
 \end{itemize}
In \cite{algintdif}, it was conjectured that $\gldim (\mI_n)=n$  and
 $\gldim (\mA_n)=n$. The aim of the paper is to prove these conjecturers.
  \begin{itemize}
\item  (Theorem \ref{AAA8Oct9}) $\gldim (\mI_n)=n$.
\item  (Theorem \ref{J18Oct9}.(2)) $\gldim (\mA_n)=n$.
 \end{itemize}
For a Noetherian algebra, its global dimension and the weak global  dimension coincide. In general, it is much more easy to compute the weak global  dimension than the global dimension as the former behaves better under many constructions (like localizations). It is a standard approach that in order to find the global dimension of a Noetherian algebra actually its weak global  dimension is computed. The algebras $\mI_n$ and $\mA_n$ are not Noetherian algebras and they contain infinite direct sums of nonzero left and right ideals. These fact are major difficulty in computing the global dimension of $\mI_n$ and $\mA_n$ as not much technique is available in the non-Noetherian case (even on the level of examples).

{\bf Ideas of the proof of Theorem \ref{AAA8Oct9}}: To show that $\gldim (\mI_n) =n$ we use an induction on $n$. The case $n=1$ is turned out to be a challenging  one. We spend entire Section \ref{GLMI1} to prove that $\gldim (\mI_1)=1$. In Section \ref{GLMI1}, we recall some facts about the algebra $\mI_1$ that are used in the proof of the case $n=1$ and in the general case. Then we use localizations and properties of certain ideals of the algebra $\mI_n$ and of the $\mI_n$-module $P_n$ (especially,  that it is a projective $\mI_n$-module) to tackle the inductive step. $\Box$

{\bf The global dimension of prime factor algebras of $\mI_n$ and $\mA_n$.} Using an explicit description of all the prime factor algebras of the algebra $\mI_n$ obtained in \cite{algintdif}  and the fact that $\gldim (\mI_n) =\gldim (\mA_n) =n$, the global dimension of prime factor algebras of the algebras $\mI_n$ and $\mA_n$ are found.
  \begin{itemize}
\item  (Theorem \ref{G18Oct9})  $\gldim (A)=n$ {\em for all prime factor algebras $A$ of} $\mI_n$.
\item  (Theorem \ref{Ja19Oct9})  $\gldim (A)=n$ {\em for all prime factor algebras $A$ of} $\mA_n$.
 \end{itemize}

 {\bf The weak global dimension of factor algebras of $\mI_n$ and $\mA_n$.} In Section \ref{GLDIMFACTOR}, a technique is developed of finding the weak global dimension of a ring via certain left/right localizations that satisfy some natural conditions (Theorem \ref{1May17} and Theorem \ref{BGanya1May17}). It is applied to  factor algebras of the algebra $\mI_n$ and $\mA_n$ to prove the following theorems.
   \begin{itemize}
\item  (Theorem \ref{29Apr17}.(2))  $\wdim (A)=n$ {\em for all  factor algebras $A$ of\,}  $\mI_n$.
\item  (Theorem \ref{An29Apr17}.(2))  $\wdim (A)=n$ {\em for all  factor algebras $A$ of} $\mA_n$.
 \end{itemize}
 Explicit descriptions of ideals of the algebras $\mA_n$ and $\mI_n$ (obtained in \cite{Bav-Jacalg,algintdif}) are one of the crucial steps in the proofs of the theorems.

{\bf An analogue of    Hilbert's
Syzygy Theorem for the algebras $\mI_n$, $\mA_n$ and their prime factor algebras.}
 Many classical algebras are tensor product of algebras (eg,
 $P_n=P_1^{\t n}$, $A_n= A_1^{\t n}$, $\mA_n= \mA_1^{\t n}$,
 $\mI_n=\mI_1^{\t n}$, etc). In general, it is difficult to find the left global
 dimension $\lgldim (A\t B)$ of the tensor product of two algebras (even to
 answer the question  when it is finite). In
 \cite{ERZ},  it was pointed out  by Eilenberg, Rosenberg and  Zelinsky that `{\em
the questions concerning the dimension of the tensor product of
two algebras have turned out to be surprisingly difficult.}' An
answer is known if one of the algebras is a polynomial algebra:
$${\bf  Hilbert's \; Syzygy\;  Theorem}:\;\;\;  \lgldim (P_n \t B) = \lgldim (P_n) +\lgldim (B) =
n+\lgldim (B). $$ In \cite{THM,glgwa}, an analogue of  Hilbert's
Syzygy Theorem was proven  for certain {\em generalized Weyl
algebras} $A$ (eg, $A=A_n$, the Weyl algebra):
$$ \lgldim (A\t B) = \lgldim (A) +\lgldim (B)$$
for all left Noetherian finitely generated algebras $B$ ($K$ is an
algebraically closed uncountable field of characteristic zero). In
this paper, the same result is proven for the algebras $\mI_n$, $\mA_n$  and
for all their prime factor algebras (Theorem \ref{G19Oct9} and Theorem \ref{J19Oct9}).

\begin{itemize}
\item  (Theorem \ref{G19Oct9}.(2)) {\em Let $K$ be an algebraically closed uncountable field of
characteristic zero and $B$ be a left Noetherian finitely
generated algebra over $K$. Then} $\lgldim (\mI_n \t B) = \lgldim (\mI_n ) +\lgldim (B)  = n+\lgldim (B)$.
\item  (Theorem \ref{J19Oct9}.(2)) {\em Let $K$ be an algebraically closed uncountable field of
characteristic zero and $B$ be a left Noetherian finitely
generated algebra over $K$. Then}   $\lgldim (\mA_n\t B) = \lgldim (\mA_n) +\lgldim (B)  = n+\lgldim (B)$.
 \end{itemize}

In case of $\wdim$, we prove even stronger results.
   \begin{itemize}
\item  (Theorem \ref{3May17}.(2)) {\em Let $K$ be an algebraically closed uncountable field of
characteristic zero and $B$ be a left Noetherian finitely
generated algebra over $K$. Then
  $\wdim (A\t B) = \wdim (A) +\wdim (B)  = n+\wdim (B)$ for all factor algebras $A$ of} $\mI_n$.
\item  (Theorem \ref{A3May17}.(2)) {\em Let $K$ be an algebraically closed uncountable field of
characteristic zero and $B$ be a left Noetherian finitely
generated algebra over $K$. Then
  $\wdim (A\t B) = \wdim (A) +\wdim (B)  = n+\wdim (B)$ for all factor algebras $A$ of} $\mA_n$.
 \end{itemize}

{\bf The global dimension and the weak global  dimension of the (largest) left  quotient ring of $\mI_1$.} Let $R$ be a ring and  $\CC_R$ be the set of regular elements of $R$ (the set of non-zero-divisors). The ring $Q_{l,cl} (R):= \CC_R^{-1}R$ (resp., $Q_{r, cl}(R)=R\CC_R^{-1}$) is called the {\em classical left} (resp., {\em right}) {\em quotient ring} of $R$. The left and right classical quotient rings not always exist. For the algebra $\mI_1$ neither the left nor right classical quotient ring exists \cite{intdifline}. For an arbitrary ring $R$, there exists the {\em largest  left} (resp., {\em right}) {\em Ore set} of $R$ that consists of regular elements of $R$, it is denoted by $S_l(R)$ (resp., $S_r(R)$) and the ring $Q_l (R):= S_l(R)^{-1}R$ (resp., $Q_r(R)=RS_r(R)^{-1}$) is called the {\em (largest)  left} (resp., {\em right}) {\em quotient ring} of $R$, \cite{intdifline,LargLQuotRing}. If the ring $Q_{l, cl}(R)$ (resp., $Q_{r, cl}(R)$) exists then $Q_{l, cl}(R)=Q_l(R)$ (resp., $Q_{r, cl}(R)=Q_r(R)$). For the algebra $\mI_1$, the sets $S_l(\mI_1)$, $S_r(\mI_1)$ and the rings $Q_l (\mI_1)$, $Q_r(\mI_1)$ are explicitly described in \cite[Theorem 9.7]{intdifline}. In particular, $S_l(\mI_1) \neq S_r(\mI_1)$ and the rings  $Q_l (\mI_1)$ and  $Q_r(\mI_1)$ are not isomorphic. For the Weyl algebra $A_1$, $Q_{l, cl}(A_1)=Q_l(A_1)$ is a skew field (a division ring), and so $\gldim (Q_l(A_1))=0$. This is not the case for the algebra $\mI_1$ as the next theorem shows.

\begin{theorem}\label{A20Apr17}\marginpar{A20Apr17}
\begin{enumerate}
\item $\lgldim (Q_l(\mI_1)) =1$.
\item For all regular left Ore sets $S$ of the algebra $\mI_1$, $\lgldim (S^{-1}\mI_1)) =1$.
\end{enumerate}
\end{theorem}

For left Noetherian rings the left global dimension coincides with the weak global  dimension. For not left Noetherian rings this is not true, in general, and the next result demonstrates this phenomenon.
\begin{theorem}\label{B20Apr17}\marginpar{B20Apr17}
$\wdim (Q_l(\mI_1)) =0$   and $\wdim (Q_r(\mI_1)) =0$.
\end{theorem}

The algebras of integro-differential operators have connections to  the
{\em Rota-Baxter} algebras. The latter appeared in the work of
Baxter \cite{Baxter} and further explored by Rota \cite{Rota-1,
Rota-2}, Cartier \cite{Cartier}, and Atkinson \cite{Atkinson}, and
more recently by many others: Aguiar, Moreira
\cite{Aguiar-Moreira}; Cassidy, Guo, Keigher, Sit, Ebrahimi-Fard
\cite{Cassidy-Guo-Keigher-Sit}, \cite{Ebrahimi-Guo}; Connes,
Kreimer, Marcoli \cite{Connes-Kreimer}, \cite{Connes-Marcoli},  Guo,
Regensburger, Rosenkranz and Middeke
 \cite{Regensburger-Rosenkranz-Middeke,Guo-Regen-Rosen-2014}, name just a  few.


\section{The global dimension of $\mI_1$ is 1}\label{GLMI1}

The aim of this section is to prove that the global dimension of
the algebra $\mI_1$ is $1$ (Theorem \ref{9Oct10}). For reader's
convenience we split the proof into several steps Lemma
\ref{a9Oct10}--Lemma \ref{g9Oct10}. Many of the steps are
interesting on their own (like Theorem \ref{A9Oct10}). We start
this section by recalling some necessary facts about the algebra
$\mI_1$, see \cite{algintdif,intdifline} for details.

The algebra $\mI_1$  is generated by the elements $\der $, $H:=
\der x$ and $\int$ (since $x=\int H$) that satisfy the defining
relations (Proposition 2.2, \cite{algintdif}): $$\der \int =
1, \;\; [H, \int ] = \int, \;\; [H, \der ] =-\der , \;\;
H(1-\int\der ) =(1-\int\der ) H = 1-\int\der .$$ The elements of
the algebra $\mI_1$,  
\begin{equation}\label{eijdef}
e_{ij}:=\int^i\der^j-\int^{i+1}\der^{j+1}, \;\; i,j\in \N ,
\end{equation}
satisfy the relations $e_{ij}e_{kl}=\d_{jk}e_{il}$ where $\d_{jk}$
is the Kronecker delta function. Notice that
$e_{ij}=\int^ie_{00}\der^j$. The matrices of the linear maps
$e_{ij}\in \End_K(K[x])$ with respect to the basis $\{ x^{[s]}:=
\frac{x^s}{s!}\}_{s\in \N}$ of the polynomial algebra $K[x]$  are
the elementary matrices, i.e.,
$ e_{ij}(x^{[s]})=\d_{js}
x^{[i]}$.
Let $E_{ij}\in \End_K(K[x])$ be the usual matrix units, i.e.,
$E_{ij}(x^s)= \d_{js}x^i$ for all $i,j,s\in \N$. Then
\begin{equation}\label{eijEij}
e_{ij}=\frac{j!}{i!}E_{ij},
\end{equation}
 $Ke_{ij}=KE_{ij}$, and
$F:=\bigoplus_{i,j\geq 0}Ke_{ij}= \bigoplus_{i,j\geq
0}KE_{ij}\simeq M_\infty (K)$, the algebra (without 1) of infinite
dimensional matrices. Using induction on $i$ and the fact that
$\int^je_{kk}\der^j=e_{k+j, k+j}$, we can easily  prove that
\begin{equation}\label{Iidi}
\int^i\der^i = 1-e_{00}-e_{11}-\cdots - e_{i-1,
i-1}=1-E_{00}-E_{11}-\cdots -E_{i-1,i-1}, \;\; i\geq 1.
\end{equation}


{\bf $\Z$-grading on the algebra $\mI_1$ and the canonical form of
an integro-differential operator}. The algebra
$\mI_1=\bigoplus_{i\in \Z} \mI_{1, i}$ is a $\Z$-graded algebra
($\mI_{1, i} \mI_{1, j}\subseteq \mI_{1, i+j}$ for all $i,j\in
\Z$) where
$$ \mI_{1, i} =\begin{cases}
D_1\int^i=\int^iD_1& \text{if } i>0,\\
D_1& \text{if }i=0,\\
\der^{|i|}D_1=D_1\der^{|i|}& \text{if }i<0,\\
\end{cases}
 $$
 the algebra $D_1:= K[H]\oplus \oplus_{i\in \N} Ke_{ii}$ is
a commutative non-Noetherian subalgebra of $\mI_1$, $ He_{ii} =
e_{ii}H= (i+1)e_{ii}$  for $i\in \N $;  $(\int^iD_1)_{D_1}\simeq
D_1$, $\int^id\mapsto d$; ${}_{D_1}(D_1\der^i) \simeq D_1$,
$d\der^i\mapsto d$,   for all $i\geq 0$ since $\der^i\int^i=1$.
 Notice that the maps $\cdot\int^i : D_1\ra D_1\int^i$, $d\mapsto
d\int^i$,  and $\der^i \cdot : D_1\ra \der^iD_1$, $d\mapsto
\der^id$, have the same kernel $\bigoplus_{j=0}^{i-1}Ke_{jj}$.
Each element $a$ of the algebra $\mI_1$ is the unique finite sum
\begin{equation}\label{acan}
a=\sum_{i>0} a_{-i}\der^i+a_0+\sum_{i>0}\int^ia_i +\sum_{i,j\in
\N} \l_{ij} e_{ij}
\end{equation}
where $a_k\in K[H]$ and $\l_{ij}\in K$. This is the {\em canonical
form} of the polynomial integro-differential operator
\cite{algintdif}.  So, the set $\{ H^j\der^i, H^j, \int^iH^j,
e_{st}\, | \, i\geq 1; j,s,t\geq 0\}$ is a $K$-basis for the
algebra $\mI_1$. The multiplication in the algebra $\mI_1$ is
given by the rule:
$$ \int H = (H-1) \int , \;\; H\der = \der (H-1), \;\; \int e_{ij}
= e_{i+1, j}, \;\; e_{ij}\int= e_{i,j-1}, \;\; \der e_{ij}=
e_{i-1, j}\;\; e_{ij} \der = \der e_{i, j+1}.$$
$$ He_{ii} = e_{ii}H= (i+1)e_{ii}, \;\; i\in \N, $$


{\bf The ideal $F$ of compact operators of $\mI_1$}.
 Let $V$ be an infinite dimensional vector space over a field $K$.
A linear map $\v \in \End_K(V)$ is called a {\em compact} linear
map/operator if $\dim_K(\im (\v ))<\infty$. The set $\CC =\CC (V)$
of all compact operators is a (two-sided) ideal of the algebra
$\End_K(V)$.
 The algebra
$\mI_1$ has the only proper ideal $F=\bigoplus_{i,j\in \N}Ke_{ij}
\simeq M_\infty (K)$, the ideal of compact operators in $\mI_1$,
$F=\mI_1\cap \CC (K[x])$,
 $F^2= F$, \cite{algintdif}. The factor algebra $\mI_1/F$ is canonically isomorphic to
the skew Laurent polynomial algebra
$$B_1:= K[H][\der, \der^{-1} ;
\tau ], \;\; \tau (H) = H+1,\;\; {\rm  via}\;\; \der \mapsto \der, \;\;  \int\mapsto
\der^{-1},\;\;  H\mapsto H$$
 where $\der^{\pm 1}\alpha = \tau^{\pm
1}(\alpha ) \der^{\pm 1}$ for all elements $\alpha \in K[H]$. The
algebra $B_1$ is canonically isomorphic to the (left and right)
localization $A_{1, \der }$ of the Weyl algebra $A_1$ at the
powers of the element $\der$ (notice that $x= \der^{-1} H$).
Therefore, they have the common skew field of fractions, $\Frac
(A_1) = \Frac (B_1)$, the {\em first Weyl skew field}. The algebra
$B_1$ is a subalgebra of the skew Laurent polynomial algebra
$$\CB_1:= K(H) [ \der, \der^{-1} ; \tau ]$$ where $K(H)$ is the
field of rational functions over the field $K$ in $H$. The algebra
$\CB_1 = S^{-1} B_1$ is the left and right localization of the
algebra $B_1$ at the multiplicative set $S=K[H]\backslash \{ 0\}$.
The algebra $\CB_1$ is a {\em noncommutative Euclidean domain},
i.e., the left and right division algorithms with remainder hold
with respect to the length function $l$ on $B_1$: $l(\alpha_m
\der^m+\alpha_{m+1}\der^{m+1}+\cdots +\alpha_n\der^n)=n-m$ where
$\alpha_i\in K(H)$, $\alpha_m\neq 0$, $\alpha_n\neq 0$, and
$m<\cdots <n$. In particular, the algebra $\CB_1$ is a principal
left and right ideal domain. A $\CB_1$-module $M$ is simple iff
$M\simeq \CB_1/\CB_1b$ for some irreducible element $b\in \CB_1$,
and $\CB_1/\CB_1b\simeq \CB_1/\CB_1c$ iff the elements $b$ and $c$
are {\em similar} (that is,  there exists an element $d\in \CB_1$
such that $1$ is the greatest common right divisor of $c$ and $d$,
and $bd$ is the least common left multiple of $c$ and $d$).


{\bf The involution $*$ on the algebra $\mI_1$}. The algebra
$\mI_1$ admits the involution $*$ over the field $K$: $\der^* =
\int$, $\int^* = \der$ and $H^* = H$, i.e., it is a $K$-algebra
{\em anti-isomorphism} ($(ab)^* = b^* a^*$) such that $a^{**} =a$.
Therefore, the algebra $\mI_1$ is {\em self-dual}, i.e., it is
isomorphic to its {\em opposite algebra} $\mI_1^{op}$. As a result, the
left and right properties of the algebra $\mI_1$ are the same.
Clearly, $e_{ij}^* = e_{ji}$ for all $i,j\in \N$, and so $F^* =
F$.


{\bf Classification of $K[H]$-torsion  simple $\mI_1$-modules}. In
\cite{intdifline}, a classification of simple $\mI_1$-modules is
given. In the proof of Theorem  \ref{9Oct10} we use only
$K[H]$-torsion simple $\mI_1$-modules. Let us recall their
classification. Since the field $K$ has characteristic zero, the
group $\langle \tau \rangle \simeq \Z$ acts freely on the set
$\Max (K[H])$ of maximal ideals of the polynomial algebra $K[H]$.
That is, for each maximal ideal $\gp \in \Max (K[H])$, its orbit
$\OO (\gp ) :=\{ \tau^i (\gp ) \, | \, i\in \Z \}$ contains
infinitely many elements. For two elements $\tau^i (\gp )$ and
$\tau^j(\gp )$ of the orbit $\OO (\gp )$ we write $\tau^i (\gp )
<\tau^j (\gp )$ if $i<j$. Let $\Max (K[H])/\langle \tau \rangle$
be the set of all $\langle \tau \rangle$-orbits in $\Max (K[H])$.
If $K$ is an algebraically closed field then $\gp = (H-\l )$, for
some scalar $\l \in K$;  $\Max (K[H])\simeq K$, $(H-\l ) \lra \l
$;  and $\Max(K[H])/\langle \tau \rangle\simeq K/ \Z$.

For any algebra $A$, we denote by $\hA$ the set of isomorphism
classes of simple $A$-modules and, for  a simple $A$-module $M$,
$[M]\in \hA$ stands for its isomorphism class. For   an
isomorphism invariant property $\CP$ of simple $A$-modules, let
$\hA (\CP ) := \{ [M]\in \hA\, | \, M$ has property $\CP\}$. The
{\em socle} $\soc (M)$ of a module $M$ is the sum of all the
simple submodules if they exist and zero otherwise. Since the
algebra $\mI_1$ contains the Weyl algebra $A_1$,  which is a
simple infinite dimensional algebra,  each nonzero  $\mI_1$-module
is necessarily an {\em infinite dimensional} module.

Since the algebra $B_1=\mI_1/F$ is a factor algebra of $\mI_1$,
there is the tautological embedding
$$ \hB1\ra \hmI1, \;\; [ M]\mapsto [M].$$
 Therefore, $\hB1\subseteq \hmI1$. In \cite[Theorem 2.1.(2)]{intdifline}, it is proven that {\em the map}
\begin{equation}\label{Max1}
\Max (K[H])/ \langle \tau \rangle\ra \hB1 (K[H]-{\rm torsion}),
\;\; [\gp ] \mapsto [ B_1/B_1\gp ],
\end{equation}
 {\em  is a bijection with
the inverse map $[N]\mapsto \supp (N) := \{
 \gp \in \Max (K[H])\, | \, \gp \cdot n_\gp =0 $ for some $0\neq n_\gp
 \in N\}$, and }
\begin{equation}\label{Max2}
  \hmI1 (K[H]-{\rm torsion})= \{ [K[x]]\}\,  \coprod  \, \hB1
(K[H]-{\rm
  torsion}).
\end{equation}


The left $\mI_1$-module $F=\bigoplus_{i\in \N} E_{\N , i}$ is a
semi-simple $\mI_1$-module where ${}_{\mI_1}E_{\N , i}:=
\mI_1e_{ii}=\bigoplus_{j\in \N} Ke_{ji}\simeq \mI_1 / \mI_1 \der
\simeq {}_{\mI_1}K[x]$ is a simple $\mI_1$-module. Let $e_n:=
e_{00}+e_{11}+\cdots + e_{n n}$ for $n\in \N$. Then $\mI_1 e_n =
\bigoplus_{i=0}^n E_{\N , i}$ is a semi-simple $\mI_1$-module of
length $n+1$. The left $\mI_1$-module $F$ has the ascending
filtration $F=\bigcup_{n\in \N} \mI_1e_n$, $ \mI_1e_0\subset
\mI_1e_1\subset \cdots \subset \mI_1e_n \subset \cdots $.

\begin{lemma}\label{a9Oct10}
For all $n\geq 1$, ${}_{\mI_1}\mI_1 = \mI_1\der^n\oplus
\mI_1e_{00}\oplus\cdots \oplus \mI_1e_{n-1, n-1}$ and
${\mI_1}_{\mI_1} = \int^n\mI_1\oplus
e_{00}\mI_1\oplus\cdots$ $ \oplus e_{n-1, n-1}\mI_1$.
\end{lemma}

{\it Proof}. It suffices to prove only that the first equality
holds since then the second is obtained from the first by applying
the involution $*$ of the algebra $\mI_1$. Notice that
$\mI_1e_{00}\oplus\cdots \oplus \mI_1e_{n-1,
n-1}=\mI_1e_{n-1}$ where $e_{n-1}:= e_{00}+e_{11}+\cdots +
e_{n-1,n-1}$ since $\mI_1e_{0i}= \mI_1 e_{ii}$. Using the equality
$\int^n\der^n= 1-e_{n-1}$ (see (\ref{Iidi})), we see that $\mI_1=
\mI_1\der^n + \mI_1e_{n-1}$.  Since $\der^n e_{n-1} =0$ and
$e_{n-1}^2= e_{n-1}$, we have $\mI_1\der^n\cap \mI_1 e_{n-1} =
(\mI_1 \der^n \cap \mI_1e_{n-1})e_{n-1}\subseteq \mI_1\der^n
e_{n-1} =0$. Therefore, ${}_{\mI_1}\mI_1 = \mI_1\der^n\oplus
\mI_1e_{n-1}$.
 $\Box $


It follows from Lemma \ref{a9Oct10} that the simple $\mI_1$-module
${}_{\mI_1}K[x]\simeq \mI_1 / \mI_1\der\simeq \mI_1e_{00}$ is a
{\em projective} $\mI_1$-module. For an algebra $A$ and its
element $a$, let $a\cdot $ and $\cdot a$  be the left and right
multiplications by the element $a$ in $A$.

\begin{lemma}\label{b9Oct10}
\begin{enumerate}
\item For all elements $a\in \mI_1$, the left ideal $\mI_1a$ of
$\mI_1$ is a projective left $\mI_1$-module.  \item For all
elements $a\in \mI_1\backslash F$, the left ideal
$\ker_{\mI_1}(\cdot a)$ of $\mI_1$  is a projective left
$\mI_1$-module.
\end{enumerate}
\end{lemma}

{\it Proof}. If $a\in F$ then $\mI_1 a \subseteq F$. Since
${}_{\mI_1}F\simeq K[x]^{(\N )}$ is a semi-simple projective
$\mI_1$-module so is $\mI_1 a$. Suppose that $a\not\in F$. In
\cite{intdifline}, it is proved that $\CK := \ker_{\mI_1} (\cdot
a)$ is a finitely generated $\mI_1$-module which is contained in
$F$, and so $\CK \subseteq \mI_1 e_{n-1}$ for some $ n\geq 1$. The
$\mI_1$-module $\mI_1e_{n-1}$ is a semi-simple, hence
$\mI_1e_{n-1} = \CK \oplus L$ for some $\mI_1$-submodule $L$ of
$\mI_1e_{n-1}$. By Lemma \ref{a9Oct10},
$$ \mI_1= \mI_1\der^n \oplus \mI_1e_{n-1}= \mI_1\der^n
\oplus \CK \oplus L.$$ Therefore, the $\mI_1$-modules $\CK$
and $\mI_1\der^n \oplus L\simeq \mI_1/\CK \simeq \mI_1 a$ are
projective. $\Box $


\begin{lemma}\label{c9Oct10}
Let $V$ be a left ideal of $\mI_1$ such that $V\subseteq F$ and
$a\in \mI_1$. Then the left ideal $V+\mI_1a$ of the algebra
$\mI_1$ is a projective $\mI_1$-module. In particular, the left
ideal  $F+\mI_1 a$ of the algebra $\mI_1$ is a projective
$\mI_1$-module.
\end{lemma}

{\it Proof}. The left $\mI_1$-module $F$ is a projective
semi-simple module, hence so is its submodule $V$. There is the
short exact sequence of $\mI_1$-modules
$$ 0\ra U:= V\cap \mI_1 a \stackrel{\v}{\ra} V\oplus
\mI_1 a\stackrel{\psi}{\ra} V+\mI_1 a \ra 0$$ where $\v (u) :=
(u,-u)$ and $\psi (v,w) := v+w$. Then $V= U\oplus W$ for some,
necessarily, projective $\mI_1$-submodule  $W$ of $V$. By Lemma
\ref{b9Oct10}, the $\mI_1$-module $V\oplus \mI_1 a$ is
projective. It follows from the explicit form of the map $\v $
that ${}_{\mI_1}V\oplus \mI_1 a =\im (\v ) \oplus W\oplus
\mI_1 a$. Therefore, the $\mI_1$-module $V+\mI_1 a \simeq
(V\oplus \mI_1 a) / \im (\v ) \simeq W\oplus \mI_1 a$ is
projective. $\Box $


\begin{lemma}\label{d9Oct10}
Let $I$ be a left ideal of the algebra $\mI_1$. Then the left
$\mI_1$-module $I$ is projective iff the left ideal $F+I$ of
$\mI_1$ is a projective $\mI_1$-module.
\end{lemma}

{\it Proof}. The short exact sequence of $\mI_1$-modules
$$0\ra I \ra F+I\ra (F+I) / I \simeq F/F\cap I \ra 0$$
splits since the $\mI_1$-module $F/F\cap I$ is projective
(${}_{\mI_1}F$ is a projective semi-simple module), and the result
follows. $\Box $


\begin{lemma}\label{e9Oct10}
Let $\alpha (H) \in K[H]$ be a nonzero polynomial. Then
$\ker_{\mI_1} (\cdot \alpha (H)) = \bigoplus_{ \{ i\in \N : \alpha
(i+1) =0\} } E_{\N , i}$ is a projective semi-simple left
$\mI_1$-module of finite length $|\{ i\in \N : \alpha (i+1)
=0\}|$.
\end{lemma}

{\it Proof}. Since $\alpha \not \in F$ and $ \mI_1 / F$ is a
domain, $\ker_{\mI_1}(\cdot \alpha )= \ker_F(\cdot \alpha )$.
Since ${}_{\mI_1}F_{\mI_1}\simeq ({}_{\mI_1}K[x])\t  ( K[\der ]
_{\mI_1})$ where $ K[\der ] _{\mI_1}:= \mI_1 / \int \mI_1 = (\int
\mI_1 \oplus e_{00} \mI_1 ) / \int \mI_1 \simeq e_{00} \mI_1 =
\bigoplus_{i\in \N } Ke_{0i}=:E_{0, \N }$ (Lemma \ref{a9Oct10})
and $ e_{0i} \cdot \alpha (H) = \alpha (i+1)$,  the result is
obvious. $\Box $


\begin{theorem}\label{A9Oct10}
Let $a\in \mI_1 \backslash F$ and $I$ be a nonzero ideal of the
algebra $\mI_1$.
\begin{enumerate}
\item ${}_{\mI_1}\mI_1\simeq \ker_{\mI_1}(\cdot a) \oplus \mI_1a$.
\item The left $\mI_1$-module $I$ is projective iff so is $Ia$.
\item Suppose that the left ideal $I\cap F$ is  finitely generated. Then
\begin{enumerate}
\item $I=I\cap F\oplus I_1$ for a left ideal  $I_1$ of $\mI_1$; and  $I_1\simeq I/I\cap F$.
\item The $\mI_1$-module $I$ is projective iff the $\mI_1$-module $I_1$ is so.
\end{enumerate}
\end{enumerate}
\end{theorem}

{\it Proof}. 1. Since $a\in \mI_1 \backslash F$, by \cite[Theorem 4.2.(1)]{intdifline}, the kernel $\CK := \ker_{\mI_1} (\cdot a)$ is a finitely
generated (necessarily, semi-simple) submodule of the semi-simple
$\mI_1$-module $F= \bigoplus_{i\in \N} E_{\N , i}$ and so $\CK
\subseteq \bigoplus_{i=0}^n E_{\N , i}$ for some $n$. By Lemma
\ref{a9Oct10}, ${}_{\mI_1}\mI_1= \CK\oplus J$ for some left
ideal $J$ of $\mI_1$. Since $J\simeq \mI_1 / \CK \simeq \mI_1 a$,
statement 1 follows.

2. The intersection $\CK'=\CK\cap I$ is a finitely generated submodule of $F$ (see the proof of statement 1). We keep the notation of the proof of statement 1. Then $\mI_1=\CK'\oplus \mJ$ for some left ideal $\mJ$ of $\mI_1$. Hence, $I=I\cap \mI_1=I\cap (\CK'\oplus \mJ )=\CK'\oplus I'$ where $I'=I\cap \mJ$. There is a short exact sequence of $\mI_1$-modules, $0\ra \CK'\ra I\stackrel{\cdot a}{\ra} Ia\ra 0$. In particular, $Ia\simeq I/\CK'\simeq I'$. So, $I\simeq \CK'\oplus Ia$, and statement 2 follows.

 3(a). Since the left ideal $I_0:=I\cap F$ is finitely generated, there is a natural number $n\geq 1$ such that $I_0\subseteq V:=\bigoplus_{i=0}^{n-1}\mI_1e_{ii}$. The module $V$ is a projective  semisimple $\mI_1$-module, hence $V=I_0\oplus U$ for some projective semisimple  $\mI_1$-module $U$.  By  Lemma \ref{a9Oct10}, $\mI_1=V\oplus \mI_1\der^n$, and so $I=I_0\oplus I_1$ where $I_1=I\cap (U\oplus \mI_1\der^n)$.

3(b). The statement (b) follows from the statement (a).  $\Box $



\begin{lemma}\label{g9Oct10}
Each simple $K[H]$-torsion $\mI_1$-module has projective dimension
$\leq 1$. In more detail, $\pd_{\mI_1}(K[x])=0$ and $\pd_{\mI_1}
(M)= 1$ for all $M\in \widehat{B}_1(K[H]-{\rm torsion})$.
\end{lemma}

{\it Proof}. By (Theorem 2.1, \cite{intdifline}), each simple
$K[H]$-torsion $\mI_1$-module $M$ is isomorphic either to
${}_{\mI_1}K[x]=\mI_1/ \mI_1 \der = (\mI_1 \der \oplus \mI_1
e_{00})/ \mI_1 \der \simeq \mI_1 e_{00}=\bigoplus_{i\in \N}
Ke_{i0}$ or to $B_1/ B_1p= \mI_1/ (F+\mI_1 p)$ for some
irreducible polynomial $p\in K[H]$. Clearly, $\pd_{\mI_1}
(K[x])=0$ and $\pd_{\mI_1} (F+\mI_1 p) = \pd_{\mI_1} (\mI_1 p ) =
\pd_{\mI_1}(\mI_1) =0$, by Lemma \ref{d9Oct10} and Theorem
 \ref{A9Oct10}.(2).  The ideal $F$ is an essential left
 $\mI_1$-submodule of ${}_{\mI_1}\mI_1$ (\cite[Lemma 5.2.(2)]{algintdif}), then so is the
 left ideal $F+\mI_1p$, i.e., the projective resolution for the
 $\mI_1$-module $\mI_1 / (F+\mI_1 p)$,
 $$ 0\ra F+\mI_1 p \ra \mI_1 \ra \mI_1 / (F+\mI_1 p)\ra 0,$$
 does not split. Then $\pd_{\mI_1} (\mI_1 / (F+\mI_1 p)) =1$.
 $\Box $


By Lemma \ref{g9Oct10}, $\lgldim (\mI_1)\geq 1$. The next theorem shows that $\lgldim (\mI_1)=1$.

\begin{theorem}\label{9Oct10}
$\gldim (\mI_1)  = 1$.
\end{theorem}

{\it Proof}. The algebra $\mI_1$ is self-dual, so it suffices to
prove that $\lgldim (\mI_1) =1$, i.e., every nonzero left ideal $I$
of the algebra $\mI_1$ is a projective $\mI_1$-module. By Lemma
\ref{d9Oct10}, we may assume that $F\subseteq I$. Since
${}_{\mI_1}F\simeq K[x]^{(\N )}$ is a projective $\mI_1$-module,
we may assume that $F\subsetneqq I$. Then its image $\pi (I)$
under the algebra epimorphism $\pi : \mI_1 \ra \mI_1 / F=B_1$, $
a\mapsto \overline{a} := a+F$ is a nonzero left ideal of the
algebra $B_1$. Since the algebra $\CB_1= S^{-1}B_1$, where $S:=
K[H]\backslash \{ 0\}$, is a principal left  ideal domain, $\CB_1
\pi (I) = \CB_1 b$ for some nonzero element $b\in K[H][\der ; \tau
]\subseteq B_1$. Notice that the algebras $B_1$ and $\CB_1$ are
left and right Noetherian algebras, $B_1b\subseteq \pi (I)$ and
$S^{-1} B_1 b= \CB_1 b = S^{-1} \pi (I)$. So, the $B_1$-module
$\pi (I) / B_1 b$ is a finitely generated $K[H]$-torsion
$B_1$-module. Therefore,
$$\pi (I) = B_1 b+\sum_{i=1}^s
B_1\alpha_i^{-1} a_i b$$ for some elements $\alpha_i\in S$ and
$a_i\in K[H] [ \der ; \tau ]\subseteq B_1$ such that
$\alpha_i^{-1} a_i b\in B_1$. Clearly, $\alpha_i^{-1} a_i b\in
K[H] [ \der ; \tau ]$. Since $\CB_1 = S^{-1} B_1 = B_1S^{-1}$,
there exists elements $c_i\in K[H][\der ; \tau ]$ and $\alpha \in
S$ such that $\alpha_i^{-1} a_i = c_i \alpha^{-1}$.
$$ I = F+\mI_1 b +\sum_{i=1}^s \mI_1 c_i \alpha^{-1} b$$
where all elements $c_i\alpha^{-1} b = \alpha_i^{-1} a_i b_i\in
K[H] [ \der ; \tau ]\subseteq \mI_1$. Fix an element $\beta \in S$
such that $\alpha^{-1} b = d\beta^{-1}$ for some $d\in K[H][\der ;
\tau ]$. Notice that $\mI_1 b = \mI_1 \alpha \alpha^{-1} b = \mI_1
\alpha d \beta^{-1}$. By Theorem \ref{A9Oct10}, ${}_{\mI_1}I$ is
projective iff ${}_{\mI_1} I\beta= F\beta +\mI_1 \alpha d
+\sum_{i=1}^s \mI_1 c_i d$ is projective iff $(\mI_1 \alpha
+\sum_{i=1}^s \mI_1 c_i)d$ is projective $\mI_1$-module (Theorem
\ref{d9Oct10}) iff $\mI_1 \alpha +\sum_{i=1}^s \mI_1 c_i$ is a
projective $\mI_1$-module (Theorem \ref{A9Oct10}.(2)) iff $J:=
F+\mI_1\alpha +\sum_{i=1}^s\mI_1 c_i$ is a projective
$\mI_1$-module (Lemma \ref{d9Oct10})  iff $\pd_{\mI_1}
(M)\leq 1$ where $M:= \mI_1 / J$. If $M=0$, i.e., $J =
\mI_1$, we are done. So, we may assume that $M\neq 0$.
The nonzero $B_1$-module $N:=
\mI_1 / (F+\mI_1 \alpha )\simeq B_1/ B_1 \alpha$ is a
$K[H]$-torsion $B_1$-module of {\em finite length} (it is well
known that any proper factor module of $B_1$ has finite length).
The module $M$ is an epimorphic image of the $B_1$-module $N$.
Therefore, ${}_{B_1}M$ is $K[H]$-torsion of finite length. By
Lemma \ref{g9Oct10}, $\pd_{\mI_1} (M)\leq 1$. Therefore, $\lgldim
(\mI_1) =1$. $\Box $

{\bf The global dimension of localizations of $\mI_1$ is 1.}

{\bf Proof of Theorem \ref{A20Apr17}.} 1. Since $\lgldim (\mI_1)=1$, we must have $\lgldim (Q_l(\mI_1))\leq \lgldim (\mI_1)=1$. By \cite[Corollary 3.3.(8)]{algintdif} and \cite[Corollary 8.5]{intdifline}, the unique proper ideal $\CI$ of the algebra $Q_l(\mI_1)$ is an essential left ideal, hence $\pd_{Q_l(\mI_1)}(Q_l(\mI_1)/\CI ) \geq 1$, and so $\lgldim (Q_l(\mI_1))\geq 1$. Therefore, $\lgldim (Q_l(\mI_1))= 1$.

2. Let $S$ be a regular left Ore set of the algebra $\mI_1$. Then $S\subseteq S_l(\mI_1)$, and so $S^{-1}\mI_1\subseteq S_l(R)^{-1}R = Q_l(R)$. Now,
$$ 1=\lgldim (Q_l(\mI_1)) \leq \lgldim (S^{-1}\mI_1) \leq \lgldim (\mI_1)=1,$$
and so $  \lgldim (S^{-1}\mI_1)=1$. $\Box$

\begin{corollary}\label{a20Apr17}
\begin{enumerate}
\item $\rgldim (Q_r(\mI_1)) =1$.
\item For all regular right Ore sets $T$ of the algebra $\mI_1$, $\rgldim (\mI_1T^{-1})) =1$.
\end{enumerate}
\end{corollary}

{\it Proof.} 1. $\rgldim (Q_r(\mI_1))=\lgldim (Q_r(\mI_1)^*)=\lgldim (Q_l(\mI_1))=1$, by Theorem \ref{A20Apr17}.(1).

2.  $\rgldim (\mI_1T^{-1})=\lgldim ((\mI_1T^{-1})^*)=\lgldim ((T^*)^{-1}\mI_1)=1$, by Theorem \ref{A20Apr17}.(2) since $T^*$ is a regular left Ore set of $\mI_1$. $\Box$

{\bf Proof of Theorem \ref{B20Apr17}.} By \cite[Corollary 4.5]{THM},
$$\wdim (Q_l(\mI_1))=\max \{ \wdim (S_\der^{-1}Q_l(\mI_1)), \nu \}\;\; {\rm  where} \;\; \nu = \sup \{ \fd_{\mI_1}(M)\, | \, S_\der^{-1}M=0\}.$$ Since the algebra $S_\der^{-1}Q_l(\mI_1)\simeq Q_{l,cl}(A_1)$ is a division ring, by \cite[Corollary 8.5.(7b]{intdifline},
its weak global  dimension is 0. The $\mI_1$-module $K[x]$ is projective, hence so is the $Q_l(\mI_1)$-module $S_l(\mI_1)^{-1}K[x]$. In particular, the flat dimension of the $Q_l(\mI_1)$-module $S_l(\mI_1)^{-1}K[x]$ is 0. Every nonzero $S_\der$-torsion $Q_l(\mI_1)$-module is a direct sum of copies of the projective $Q_l(\mI_1)$-module $S_l(\mI_1)^{-1}K[x]$. Hence, $\nu =0$, and so $\wdim (Q_l(\mI_1))=\max \{ 0,0\}=0$.

Now, $\wdim (Q_r(\mI_1))=\wdim (Q_l(\mI_1)^*)=\wdim (Q_l(\mI_1))=0.$  $\Box$


\section{The  global dimension of the algebra
$\mI_n$}\label{WGDAI}

In this section, we prove that the global dimension
of the algebra $\mI_n$ and of all its prime factor algebras is
$n$ (Theorem \ref{AAA8Oct9}). An analogue of  Hilbert's Syzygy
Theorem  is proven for the algebra $\mI_n$ and for all its
prime factor algebras (Theorem \ref{G19Oct9}).

{\bf Classification of ideals of the algebra $\mI_n$.}  Theorem \ref{b10Oct9} describes all the ideals of the algebra $\mI_n$ and  shows that the ideal
theory of $\mI_n$ is `very arithmetic' (it is the best possible and `simple' ideal theory one can imagine).  Let $\CB_n$ be the set of
all functions $f:\{ 1, 2, \ldots , n\} \ra  \{ 0,1\}$. For each
function $f\in \CB_n$, $I_f:= I_{f(1)}\t \cdots \t I_{f(n)}$ is
the ideal of $\mI_n$ where $I_0:=F$ and $I_1:= \mI_1$.  Let
$\CC_n$ be the set of all  subsets of $\CB_n$ all distinct
elements of which are incomparable (two distinct elements $f$ and
$g$ of $\CB_n$ are {\em incomparable} if neither $f(i)\leq  g(i)$
nor $f(i)\geq g(i)$ for all $i$). For each $C\in \CC_n$, let
$I_C:= \sum_{f\in C}I_f$, the ideal of $\mI_n$. The number $\gd_n$
of elements in the set $\CC_n$ is called the {\em Dedekind
number}. It appeared in the paper of Dedekind
\cite{Dedekind-1871}. An asymptotic of the Dedekind numbers was
found by Korshunov \cite{Korshunov-1977}.

\begin{theorem}\label{b10Oct9}
{\rm (\cite[Corollary 3.3]{algintdif}.)}
\begin{enumerate}
\item The algebra $\mI_n$ is a prime algebra. \item The set of
height one prime ideals of the algebra $\mI_n$ is $\{ \gp_1:=
F\t\mI_{n-1}, \gp_1:= \mI_1\t F\t\mI_{n-2},\ldots , \gp_n:=
\mI_{n-1}\t F\}$. \item Each ideal of the algebra $\mI_n$ is an
idempotent ideal ($\ga^2= \ga$). \item The ideals of the algebra
$\mI_n$ commute ($\ga \gb = \gb \ga$). \item The lattice $\CJ
(\mI_n)$ of ideals of the algebra $\mI_n$ is distributive. \item
The classical Krull dimension $\clKdim (\mI_n)$ of the algebra
$\mI_n$ is $n$. \item $\ga \gb = \ga \cap \gb$ for all ideals $\ga
$ and $\gb $ of the algebra $\mI_n$. \item The ideal $\ga_n :=
\gp_1+\cdots + \gp_n$ is the largest (hence, the only maximal)
ideal of $\mI_n$ distinct from $\mI_n$, and $F_n = F^{\t
n}=\bigcap_{i=1}^n \gp_i$ is the smallest nonzero ideal of
$\mI_n$. \item {\rm (A classification of  ideals of $\mI_n$)} The
map
 $\CC_n\ra \CJ (\mI_n)$, $C\mapsto I_C:= \sum_{f\in C}I_f$
 is a bijection where $I_\emptyset :=0$. The number of ideals of
 $\mI_n$ is the Dedekind number $\gd_n$.  Moreover,
$2-n+\sum_{i=1}^n2^{n\choose i}\leq \gd_n \leq 2^{2^n}$. For
$n=1$, $F$ is the unique proper ideal of the algebra $\mI_1$.
\item {\rm (A classification of prime ideals of $\mI_n$)} Let
$\Sub_n$ be the set of all subsets of $\{ 1, \ldots , n\}$. The
map $\Sub_n\ra \Spec (\mI_n)$, $ I\mapsto \gp_I:= \sum_{i\in
I}\gp_i$, $\emptyset \mapsto 0$, is a bijection, i.e., any nonzero
prime ideal of $\mI_n$ is a unique sum of primes of height 1;
$|\Spec (\mI_n)|=2^n$; the height of $\gp_I$ is $| I|$; and
 \item  $\gp_I\subset \gp_J$ iff $I\subset
 J$.
\end{enumerate}
\end{theorem}

{\bf The  weak dimension of the algebra $\mS_n$}. Let $S$ be a
non-empty multiplicatively closed subset of a ring $R$, and let
$\ass (S):= \{ r\in R\, | \, sr =0$ for some $s\in S\}$. Then a
{\em left quotient ring} of $R$ with respect to $S$ is a ring $Q$
together with a homomorphism $\v : R\ra Q$ such that

(i) for all $s\in S$, $\v (s)$ is a unit in $Q$,

(ii) for all $q\in Q$, $q= \v (s)^{-1} \v (r)$ for some $r\in R$
and $s\in S$, and

(iii) $\ker (\v ) = \ass (S)$.

If there exists a left quotient ring $Q$ of $R$ with respect to
$S$ then it is unique up to isomorphism, and it is denoted by
$S^{-1}R$. It is also said that the ring $Q$ is the {\em left
localization} of the ring $R$ at $S$. The ring $S^{-1}R$ exists iff $S$ is a {\em left denominator set} of the ring $R$, that is $S$ is a {\em left Ore set} such that $rs=0$ for some elements $r\in R$ and $s\in S$ implies $s'r=0$ for some element $s'\in S$ ($S$ is a left Ore set if $Sr\cap Rs\neq \emptyset$ for all elements $r\in R$ and $s\in S$). If $M$ is a left $R$-module then  the set $\tor_S(M):=\{ m\in M\, | \, sm=0$ for some $s\in S\}$ is a submodule of $M$ which is called the $S$-{\em torsion submodule} of $M$. For a left denominator set $S$ of $R$, $\ass (S) = \ass_{l,R}(S):=\{ r\in R\, | \, sr=0$ for some element $s\in S\}$ is an ideal. In a similar way, a right denominator set of $R$ is defined. If $T$ is a right denominator set of $R$ then the ring $RT^{-1} = \{ rt^{-1}\, | \, r\in R, t\in T\}$ is a localization of $R$ at $T$ on the right and the set $\ass_{r, R}(T):=\{ r\in R\, | \, rt=0$ for some element $t\in T\}$ is an ideal of $R$.

{\it Example 1}. Let $S:= S_\der:= \{ \der^i, i\geq 0\}$ and
$R=\mI_1$. Then  $\ass (S) =F$, $\mI_1/\ass (S)= B_1$ and the
conditions (i)-(iii) hold where $Q=B_1$. This means that the ring
$B_1=\mI_1/F$ is the left quotient ring of $\mI_1$ at $S$,
 $B_1\simeq S^{-1}_\der \mI_1$.

{\it Example 2}. Let $S:=S_{\der_1, \ldots , \der_n}:= \{
\der^\alpha , \alpha \in \N^n\}$ and $R=\mI_n$. Then 
\begin{equation}\label{Sy1yn}
S_{\der_1, \ldots , \der_n}^{-1}\mI_n\simeq \bigotimes_{i=1}^n S_{\der_i}^{-1}\mI_1(i)=\bigotimes_{i=1}^n B_1(i)=:B_n
\end{equation}
where $B_i=K[H_i][\der_i, \der_i^{-1}; \tau_i]$ is a skew Laurent polynomial ring where $\tau_i(H_i) = H_i+1$.
Furthermore, $\ass
(S_{\der_1, \ldots , \der_n}) = \ga_n$ and   $\mI_n / \ga_n = B_n$, see \cite{algintdif} for details (where $\ga_n$ is the unique maximal ideal of the algebra $\mI_n$). The right localization $\mI_n
S_{\der_1, \ldots , \der_n}^{-1}$ of $\mI_n$ at $S_{\der_1, \ldots
, \der_n}$ does not exist, \cite{algintdif}. By
applying the involution $*$ to  (\ref{Sy1yn}), we see that
\begin{equation}\label{1Sy1yn}
\mI_nS^{-1}_{\int_1, \ldots , \int_n}\simeq B_n.
\end{equation}
So,  the algebra $B_n$ is the right localization of $\mI_n$ at the
multiplicatively closed set $S_{\int_1, \ldots , \int_n}:= \{
\int^\alpha:= \int_1^{\alpha_1}\cdots \int_n^{\alpha_n} \, | \, \alpha =(\alpha_1, \ldots , \alpha_n)\in \N^n\}$.

$\noindent $

Given  a ring $R$ and modules ${}_RM$ and $N_R$, we denote by $\pd
({}_RM)$ and $\pd (N_R)$ their projective dimensions. Let us
recall two results which will be used  in  proofs later.

\begin{proposition}\label{Aus-sub}
(\cite{Aus-NagoyaMJ-55}.) Let $M$ be a module over an algebra
$A$, $I$ a non-empty well-ordered set, $\{ M_i\}_{i\in I}$ be  a
family of submodules of $M$ such that if $i,j\in I$ and $i\leq j$
then $M_i\subseteq M_j$. If $M=\bigcup_{i\in I}M_i$ and
$\pd_A(M_i/M_{<i})\leq n$ for all $i\in I$ where $M_{<i}:=
\bigcup_{j<i}M_j$ then $\pd_A(M)\leq n$.
\end{proposition}
Let $V\subseteq U\subseteq W$ be modules. Then the factor module
$U/V$ is called a {\em sub-factor} of the module $W$. Let $\wdim$
and $\fd$ denote  the  {\em  weak global dimension} and the {\em
flat} dimension,  respectively.

It is obvious that $P_n\simeq A_n/\sum_{i=1}^n A_n \der_i$ where $A_n$ is the Weyl algebra. A
similar result is true for the $\mI_n$-module $P_n$ (Proposition
\ref{C17Oct9}.(2)). Notice that $\pd_{A_n}(P_n)=n$ but
$\pd_{\mI_n}(P_n)=0$ (Proposition \ref{C17Oct9}.(3)).

\begin{proposition}\label{C17Oct9}
(\cite[Proposition 6.1]{algintdif}.)\begin{enumerate}
\item $\mI_1=\mI_1\der \bigoplus \mI_1e_{00}$ and $\mI_1= \int
\mI_1\bigoplus e_{00}\mI_1$. \item ${}_{\mI_n}P_n\simeq \mI_n /
\sum_{i=1}^n \mI_n \der_i$. \item The $\mI_n$-module $P_n$ is
projective. \item $F_n =F^{\t n}$ is a left and right projective
$\mI_n$-module. \item The projective dimension of the left and
right $\mI_n$-module $\mI_n / F_n$ is $1$. \item For each element
$\alpha \in \N^n$, the $\mI_n$-module $\mI_n / \mI_n \der^\alpha$
is projective.
\end{enumerate}
\end{proposition}

Let $A$ be an algebra, $M$ be an $A$-module and $S$ be a left denominator set of $A$. The set of $S$-torsion elements in $M$, ${\rm tor}_S(M):= \{ m\in M\, | \, sm$ for some $s\in S\}$, is an $A$-submodule of $M$.

\begin{theorem}\label{AAA8Oct9}
$\gldim (\mI_n)=n$.
\end{theorem}

{\it Proof}. The algebra $\mI_n$ is self-dual. Therefore, $\lgldim (\mI_n)=\rgldim (\mI_n)$. To prove that $\lgldim (\mI_n)=n$  we use induction on $n$. The case $n=1$ holds (Theorem \ref{9Oct10}). So, we assume that $n>1$ and the equality is true for all natural numbers $n'<n$. By \cite[Proposition 6.7]{algintdif},
$$ n\leq \lgldim (\mI_n)\leq 2n.$$
Let $I$ be a left ideal of the algebra $\mI_n$. We have to show that
$$\pdIn (I)\leq n-1$$ (since then $\lgldim (\mI_n)\leq n$, and so $\lgldim (\mI_n)=n$). The algebra $\mI_n$ is a tensor product  of algebras $\mI_1\t \mI_{n-1}$. By applying the exact functor  $-\t \mI_{n-1}$ to the short exact sequence of $\mI_1$-modules $0\ra F\ra \mI_1\ra B_1\ra 0$ we obtain the short exact sequence of $\mI_n$-modules
 $0\ra F\t \mI_{n-1}\ra \mI_n\ra B_1\t \mI_{n-1}\ra 0$. As a result we have the short exact sequence of left $\mI_n$-modules
\begin{equation}\label{1}
 0\ra I_0:=I\cap F\t \mI_{n-1}\ra I\ra \bI :=I/I_0\ra 0.
\end{equation}
 Notice that ${\rm tor}_{S_{\der_1}}(I)=I_0$ (where  $S_{\der_1} :=\{\der_1^i\, | \, i\in \N \}$) and $S_{\der_1}^{-1} I = S_{\der_1}^{-1}\bI $. Therefore,
\begin{equation}\label{2}
 0\ra \bI \ra B_1\t \mI_{n-1}\ra J:= B_1\t \mI_{n-1}/\bI \ra 0
\end{equation}
is an exact sequence of $\mI_n$-modules since $B_1\simeq S_{\der_1}^{-1}\mI_1$.

{\em Claim.} $\pdIn (I_0)\leq n-1$ {\em and} $\pdIn (J)\leq n-1$: The $\mI_n$-modules $I_0$ and $J$ are $S_{\der_1}$-torsion. By \cite[Corollary 6.3]{algintdif}, each $S_{\der_1}$-torsion $\mI_n$-module $M$ admits a family $\{ T_\l\}_{\l \in \L}$ of $ \mI_n$-modules such that $M=\cup_{\l \in \L}T_\l $ where $(\L , \leq )$ is a  well-ordered set such that $\l\leq \mu$ implies $T_\l \subseteq T_\mu$ and $T_\l /\cup_{\mu <\l } T_\mu \simeq K[x_1]\t \CT_\l$ for some $\mI_{n-1}$-module $\CT_\l$. The $\mI_1$-module $K[x_1]$ is projective (Proposition \ref{C17Oct9}.(3)). By Proposition \ref{Aus-sub},
$$\pdIn (M)\leq \max \{ \pdIn (K[x]\t \CT_\l\}_{\l \in \L}\leq \lgldim (\mI_{n-1})=n-1,$$ by induction. The proof of the Claim is complete.

Suppose that $J=0$. Then,  by (\ref{1}) and (\ref{2}),
$$ \pdIn (I)\leq \max \{ \pdIn (I_0) , \pdIn (B_1\t\mI_{n-1} )\}\leq \max \{ n-1, 1\} = n-1$$
since $\pdIn (I_0)\leq n-1$ (by the Claim) and $ \pdIn (B_1\t\mI_{n-1})\leq \pd_{\mI_1}(B_1)\leq \lgldim (\mI_1)=1$.

Suppose that $\bI =0$. Then,  by (\ref{1}),
$ \pdIn (I)=\pdIn (I_0)\leq  n-1$, by the Claim.

So, we may assume that $\bI \neq 0$ and  $J\neq 0$. Then the $\mI_n$-module $\bI$ is an essential submodule of $B_1\t \mI_{n-1}$, and so the short exact sequence (\ref{2}) does not split since $J\neq 0$. In particular, $\pdIn (J)\geq 1$.
 Since $\lgldim (\mI_1)=1$, $\pdIn (B_1\t \mI_{n-1})\leq \pd_{\mI_1}(B_1)\leq 1$. Therefore,
 $$\pdIn (\bI )\leq \max \{ \pd_{\mI_1}(B_1\t \mI_{n-1}), \pd_{\mI_1}(J)\} \leq  \max \{ 1, \pd_{\mI_1}(J)\}\leq \pdIn (J).$$
   By (\ref{1}) and (\ref{2}) and the Claim,
\begin{equation}\label{3}
\pdIn (I)\leq \max \{ \pdIn (I_0) , \pdIn (\bI )\}\leq \max \{ \pdIn (I_0) , \pdIn (J )\}\leq n-1,
\end{equation}
as required. $\Box$

\begin{theorem}\label{G18Oct9}
Let $\mI_{n,m}:= B_{n-m}\t \mI_m$ where $m=0, 1, \ldots , n$ and
$\mI_0 = B_0:=K$. Then
\begin{enumerate}
\item $\gldim (\mI_{n,m}) = n$ for all
$m=0, 1, \ldots , n$.
\item  For all prime factor algebras $A$ of $\mI_n$,  $\gldim (A) =
 n$.
\end{enumerate}
\end{theorem}

{\it Proof}. 1. The algebras $B_{n-m}$ and $\mI_n$ are self-dual,
hence so is their tensor product $\mI_{n,m}$. So, the left and right global dimensions of the algebra $\mI_{n,m}$ coincide and they are denoted by $\gldim (\mI_{n,m})$.  Recall that $\mI_{n,m}=S_{\der_1, \ldots , \der_{n-m}}^{-1}\mI_n$ and $\gldim (B_k)=k$ for all $k$. Hence,
\begin{eqnarray*}
 n&=&n-m+m=\gldim (B_{n-m})+\gldim (\mI_m)\leq \gldim (B_{n-m}\t \mI_m)\\
 &\leq & \lgldim (S_{\der_1, \ldots , \der_{n-m}}^{-1}\mI_{n,m}) \leq \gldim (\mI_{n,m})=n.
\end{eqnarray*}

2. By \cite[Corollary 3.3.(10)]{algintdif}, the algebra $A$ is isomorphic to $\mI_{n,m}$, and  statement 2 follows from statement 1. $\Box$

The next theorem is an analogue of  Hilbert's Syzygy Theorem for
the algebra $\mI_n$ and its prime factor algebras.

\begin{theorem}\label{G19Oct9}
Let $K$ be an algebraically closed uncountable field of
characteristic zero. Let $A$ be a prime factor algebra of $\mI_n$
(i.e., $A\simeq \mI_{n,m}$) and $B$ be a left Noetherian finitely
generated algebra over $K$. Then
\begin{enumerate}
\item $\lgldim (A\t B) = \lgldim (A) +\lgldim (B)  = n+\lgldim
(B)$. \item  $\lgldim (\mI_n \t B) = \lgldim (\mI_n ) +\lgldim (B)  = n+\lgldim
(B)$.
\end{enumerate}
\end{theorem}

{\it Proof}. 1. Recall that $A\simeq \mI_{n,m}$ for some $m\in \{0,
1, \ldots , n\}$ and $\lgldim (\mI_{n,m})=  n$
(Theorem \ref{AAA8Oct9}). The case $m=0$, i.e.,  $A=B_n$, is well-known \cite[Corollary 1.5]{THM}.
Suppose that $m>0$ and we assume that the result is true for all $m'<m$. Let $l=\lgldim (B)$.
 By Theorem  \ref{G18Oct9},
$$ 1\leq n+l=\lgldim (\mI_{n,m})+\lgldim (B)  \leq \lgldim
(\mI_{n,m}\t B). $$
  It suffices to
show that $\lgldim (\mI_{n,m}\t B) \leq n+\lgldim (B)$ or, equivalently, that
$\pd_C(I)\leq n+l-1$ for all   left ideals $I$  of the algebra $C:=\mI_{n,m} \t B$.
 Using the localization at $S_{\der_1}$ and repeating the arguments of  the proof of Theorem \ref{AAA8Oct9} we obtain the  inequality
 $$ \pd_C(I)\leq \max \{ \pd_C(I_0), \pd_C(J)\}$$ where $I_0={\rm tor}_{S_{\der_1}}(I)$, $J=S_{\der_1}^{-1}C/\bI$  and $\bI = C/I_0$. The $C$-modules $I_0$ and $J$ are $S_{\der_1}$-torsion, hence their projective dimensions are $\leq \lgldim (S_{\der_1}^{-1}C)-1=
\lgldim (\mI_{n, m-1})-1= n+l-1$ (by repeating the proof of the Claim of Theorem \ref{AAA8Oct9}). Hence,  $\pd_C(I)\leq n+l-1$.

2. Statement 2 is a particular case of statement 1 since $\mI_n = \mI_{n,n}$.  $\Box$


\section{The  global dimension of the Jacobian algebra
$\mA_n$}\label{WGDJJ}

In this section, we prove that the global dimension
of the Jacobian algebra $\mA_n$ and of all its prime factor
algebras is $n$ (Theorem \ref{J18Oct9}, Corollary
\ref{Ja19Oct9}). An analogue of  Hilbert's Syzygy Theorem is
 proven for the Jacobian algebras $\mA_n$ and for all its
prime factor algebras (Theorem \ref{J19Oct9}).

{\bf The involution $\th$ on $\mA_n$, \cite{Bav-Jacalg}}. The Weyl algebra $A_n$ admits the {\em involution}
$$\th : A_n\ra A_n, \;\; x_i\mapsto \der_i, \;\; \der_i\mapsto
x_i, \;\; {\rm for}\;\;  i=1, \ldots , n,$$ i.e., it is a $K$-algebra
anti-isomorphism ($\th (ab) = \th (b) \th (a)$) such that  $\th^2=
{\rm id}_{A_n}$. The involution $\th$ can be uniquely extended to
the involution of $\mA_n$ by the rule
\begin{equation}\label{thinv}
\th : \mA_n\ra \mA_n, \;\; x_i\mapsto \der_i, \;\; \der_i\mapsto
x_i,\;\;  \th (H_i^{-1})= H_i^{-1},  \;\; i=1, \ldots , n.
\end{equation}

{\bf The algebras $\CA_n$.} The  polynomial algebra  $P_n=K[H_1, \ldots , H_n]$ admits a set  $\{ \sigma_1,...,\sigma_n\} $ of  commuting automorphisms where
 $\s_i(H_i)=H_i-1$ and $\s_i(H_j)=H_j$, for $i\neq j$.  The multiplicative submonoid  $\CS_n$
of $P_n$ generated by the elements $H_i+j$ (where $i=1, \ldots , n$
and $j\in \mathbb{Z}$) is  an Ore set of the Weyl algebra $A_n$,
and the (two-sided)
localization $\CA_n:=\CS_n^{-1}A_n\simeq A_n\CS_n^{-1}$ of the Weyl algebra $A_n$ at
$\CS_n$ is the {\em skew Laurent polynomial ring}
\begin{equation}\label{Anskewlaurent}
\CA_n=\CS_n^{-1}P_n [x_1^{\pm 1}, \ldots ,x_n^{\pm 1};
\sigma_1,...,\sigma_n]
\end{equation}
with coefficients from the localization of $P_n$ at $\CS_n$,
$$ \CS_n^{-1}P_n=K[H_1^{\pm 1}, (H_1 \pm 1)^{-1}, (H_1 \pm 2)^{-1}, \ldots ,
H_n^{\pm 1}, (H_n \pm 1)^{-1}, (H_n \pm 2)^{-1}, \ldots ].$$
 We identify the Weyl algebra
$A_n$ with the subalgebra of $\CA_n$ via the monomorphism,
$$A_n\ra \CA_n, \;\; x_i\mapsto x_i,\;\; \der_i\mapsto  H_ix_i^{-1}, \;\;
i=1, \ldots , n.$$ Let $k_n=Q_{l,cl}(A_n)=Q_{r,cl}(A_n)$ be the $n$'th {\em Weyl skew field} (it exists by  Goldie's Theorem  since the Weyl algebra $A_n$ is a Noetherian
domain). Then the algebra $\CA_n$ is a $K$-subalgebra of $k_n$
generated by the elements $x_i$, $x_i^{-1}$, $H_i$ and $H_i^{-1}$,
$i=1, \ldots , n$ since, for all natural numbers $j$,
$$ (H_i\mp j)^{-1}=x_i^{\pm j}H_i^{-1}x_i^{\mp j}, \;\; i=1, \ldots , n.$$
Clearly, $\CA_n \simeq \CA_1 \t \cdots \t \CA_1$ ($n$ times).

The algebra $\mA_n$ contains a unique maximal ideal $\ga_n$ (see \cite[Corollary  2.7.(4)]{Bav-Jacalg}) and   $\mA_n/\ga_n \simeq \CA_n$ (see \cite[Eq. (22)]{Bav-Jacalg}).  Since $\th (\ga_n) = \ga_n$, the algebra $\CA_n$ inherits the involution $\th$ of the algebra $\mA_n$ which is given by the rule $\th (a+\ga_n)= \th (a) +\ga_n$ for all elements $a\in \mA_n$. In particular, the algebras $\mA_n$ and $\CA_n$ are self-dual.

 A $K$-algebra $R$ has the {\em endomorphism
property over} $K$ if, for each simple $R$-module $M$, $\End_R(M)$
is algebraic over $K$.

\begin{theorem}\label{bigAn}
{\rm (\cite{Bav-Ann2001}.)}  Let $K$ be a field of characteristic
zero.
\begin{enumerate}
\item The algebra $\CA_n$ is a simple, affine, Noetherian domain.
\item The Gelfand-Kirillov dimension $\GK (\CA_n)=3n$ $(\neq
2n=\GK (A_n))$. \item The (left and right) global dimension ${\rm
gl.dim} (\CA_n)=n$. \item The (left and right) Krull dimension
${\rm K.dim} (\CA_n)=n$. \item Let ${\rm d}={\rm gl.dim}$ or ${\rm
d}= {\rm K.dim}$. Let $R$  be a Noetherian $K$-algebra with ${\rm
d}(R)<\infty $ such that $R[t]$, the polynomial ring in a central
 indeterminate, has the endomorphism property over $K$. Then
${\rm d}(\CA_1\t R)= {\rm d}(R)+1$. If, in addition, the field $K$
is algebraically closed and uncountable, and the algebra $R$ is
affine, then
 ${\rm d}(\CA_n\t R)= {\rm d}(R)+n$.
\end{enumerate}
\end{theorem}

$\GK ({\cal A}_1)=3$ is due to A. Joseph \cite{Jos1}, p. 336;
 see also  \cite{KL}, Example 4.11, p. 45.

{\bf The Jacobian algebra $\mA_n$ is a localization of the algebra
$\mI_n$}. The algebra $\mA_n$ is a (left and right) localization of the algebra $\mI_n$,  \cite[Eq. (31)]{algintdif},    
\begin{equation}\label{mAnlIn}
\mA_n=S^{-1}\mI_n = \mI_n S^{-1},
\end{equation}
 at the (left and right) Ore set  $S:=
\{ \prod_{i=1}^n(H_i+\alpha_i)_*^{n_i} \, | \, (\alpha_i) \in
\Z^n, (n_i)\in \N^n \}$ of $\mI_n$ that consists of regular elements of $\mI_n$ where
$$(H_i+\alpha_i)_*:=\begin{cases}
H_i+\alpha_i & \text{if }\alpha_i\geq 0,\\
(H_i+\alpha_i)_1& \text{if }\alpha_i<0,\\
\end{cases}\;\; {\rm and}\;\; (H_i+\alpha_i)_1:= H_i+\alpha_i-1+x_iH_i^{-1}\der_i. $$

The left (resp. right)  localization of the Jacobian algebra
$$\mA_n = K\langle y_1, \ldots , y_n, H_1^{\pm 1}, \ldots ,
H_n^{\pm 1}, x_1, \ldots , x_n\rangle  \;\; ({\rm where}\;\; y_i:=
H_i^{-1}x_i) $$ at the left denominator  set $S_{y_1, \ldots
, y_n}:=\{ y^\alpha \, | \, \alpha \in \N^n\}$ (resp., the right denominator set  $S_{x_1,
\ldots , x_n}:=\{ x^\alpha \, | \, \alpha \in \N^n\}$) is the
algebra, \cite[Eq. (33)]{algintdif},
\begin{equation}\label{mAny}
\CA_n\simeq S_{y_1, \ldots , y_n}^{-1}\mA_n\simeq \mA_n S_{x_1,
\ldots , x_n}^{-1}.
\end{equation}

The algebras $\mA_n$ and $\CA_n$ are self-dual, hence so are
the algebras $\mA_{n,m}:=
\CA_{n-m}\t \mA_m$ where $m=0, 1, \ldots , n$ and
$\mA_0=\CA_0:=K$, and so $\lgldim (\mA_{n,m})= \rgldim (\mA_{n,m}):=
\gldim (\mA_{n,m})$.

\begin{theorem}\label{J18Oct9}
Let $\mA_{n,m}:= \CA_{n-m}\t \mA_m$ where $m=0, 1, \ldots , n$ and
$\mA_0 = \CA_0:=K$. Then
\begin{enumerate}
\item $\gldim (\mA_{n,m}) =  n$ for all
$m=0, 1, \ldots , n$.
\item  $\gldim (\mA_n) =  n$.
\end{enumerate}
\end{theorem}

{\it Proof}. 1.Recall that  $\lgldim (\mA_{n,m}) = \rgldim (\mA_{n,m}) $.
 By Theorem \ref{bigAn}.(1,3),
\begin{eqnarray*}
 n &=& \lgldim (\CA_n) =\lgldim (S^{-1}_{y_1, \ldots ,
y_n} \mA_n)\leq \lgldim (\mA_{n,m})\leq \lgldim (\mA_n)\\
  &=&  \lgldim (S^{-1} \mI_n) \leq \lgldim (\mI_n)=n\;\;\; ({\rm Theorem} \;
\ref{AAA8Oct9}).
\end{eqnarray*}
Therefore, $\lgldim (\mA_{n,m})=n$ for all $n$
and $m$.

2. Statement 2 is a particular case of statement 1 when $n=m$. $\Box $

\begin{corollary}\label{Ja19Oct9}
Let $A$ be a prime factor algebra of the algebra $\mA_n$. Then
 $\lgldim (A) = \rgldim (A) =n$.
\end{corollary}

{\it Proof}. By  \cite[Corollary 3.5]{Bav-Jacalg},  the algebra $A$
is isomorphic to the algebra $\mA_{n,m}$ for some $m$. Now, the
corollary follows from Theorem \ref{J18Oct9}. $\Box $

The next theorem is an analogue of  Hilbert's Syzygy Theorem for
the Jacobian algebras and their prime factor algebras.

\begin{theorem}\label{J19Oct9}
Let $K$ be an algebraically closed uncountable field of
characteristic zero. Let $A$ be a prime factor algebra of $\mA_n$ and $B$ be a left Noetherian finitely
generated algebra over $K$. Then
\begin{enumerate}
\item  $\lgldim (A\t B) = \lgldim (A) +\lgldim (B)  = n+\lgldim
(B)$.
\item  $\lgldim (\mA_n\t B) = \lgldim (\mA_n) +\lgldim (B)  = n+\lgldim (B)$.
\end{enumerate}
\end{theorem}

{\it Proof}. 1. Recall that $A\simeq \mA_{n,m}$ for some $m\in \{0,
1, \ldots , n\}$ and $\lgldim (\mA_{n,m})= n$
(Theorem \ref{J18Oct9}). Since
\begin{eqnarray*}
n+\lgldim (B) &=&\lgldim (\mA_{n,m})+\lgldim(B) \leq
\lgldim (\mA_{n,m}\t B) \\
&\leq & \lgldim(S^{-1}\mI_n\t B)\\
&\leq & \lgldim (\mI_n \t B) = n+\lgldim (B)\;\; ({\rm by \;
Theorem\; \ref{G19Oct9}}).
\end{eqnarray*}
Therefore, $ \lgldim (\mA_{n,m}\t
B)=n+\lgldim (B)$. The proof of statement 1 is complete.

2. Statement 2 is a particular case of statement 1 when $A=\mA_n$.  $\Box $


\section{The weak  global dimension of factor algebras of $\mI_n$}\label{GLDIMFACTOR}

The aim of this section is to show that the weak global dimension of all factor algebras of the algebra $\mI_n$ is $n$ (Theorem \ref{29Apr17}) and that an analogue of Hilbert's Syzygy Theorem holds for them (Theorem \ref{3May17}).

Let $I$ be an ideal of the algebra $\mI_n$ and $A=\mI_n/I$. The ideal $I$ is the intersection of the minimal primes over the ideal $I$,

\begin{equation}\label{IIq}
I=\bigcap_{\gq \in \min (I)}\gq  \;\;\; ({\rm see\;\;  \cite[Corollary \; 3.4.(4)]{algintdif}}).
\end{equation}

The map

\begin{equation}\label{IIq1}
f: A= \mI_n/I\ra \widetilde{A} := \prod_{\gq \in \min (I)}\mI_n/ \gq ,  \;\; a\mapsto (\ldots , a+\gq , \ldots )_{\gq \in \min (I)}
\end{equation}
is an algebra monomorphism.

For a ring $R$, an element $r\in R$ is called a {\em left} (resp., {\em right}) {\em regular} if the map $\cdot r  : R\ra R$, $ s\mapsto sr$ (resp., $r\cdot   : R\ra R$, $ s\mapsto rs$) is an injection. The set  of all left (resp., right) regular elements of the ring $R$ is denoted by ${}'\CC_R$ (resp., $\CC_R'$). An element $r\in R$ is called a {\em regular element} if it is  left and right regular. By the very definition, the set $\CC_R$ of regular elements of $R$ is equal to the intersection ${}'\CC_R\cap\CC_R'$.

\begin{proposition}\label{A30Apr17}
Let $I$ be an ideal of the algebra $\mI_n$ and $A= \mI_n / I$.
\begin{enumerate}
\item The involution $*$ on the algebra $\mI_n$ induces the involution $*$ on the algebra $A$ by the rule $(a+I)^*= a^*+I$ for all elements $a\in \mI_n$.
\item $\der_1, \ldots , \der_n\in {}'\CC_A$.
\item $\int_1, \ldots , \int_n\in \CC_A'$.
\item Let $S$ (resp., $T$) be the multiplicative monoid generated by the elements $\der_{i_1}, \ldots , \der_{i_m}$ (resp., $\int_{i_1}, \ldots , \int_{i_m}$). Then $S$ (resp., $T$) is a left (resp., right) denominator set of the algebra $A$,  $\ass_{l, A}(S)=\ass_{r, A}(T)=(\gp_{i_1}+\cdots + \gp_{i_m}+I)/I$ and $S^{-1}A\simeq \mI_n / (\gp_{i_1}+\cdots + \gp_{i_m}+I)\simeq AT^{-1}$.
\end{enumerate}
\end{proposition}

{\it Proof}. 1. By \cite[Lemma 5.1.(1)]{algintdif}, all the ideals of the algebra $\mI_n$ are $*$-invariant, and so statement 1 follows.

2. For every prime ideal $\gp$ of the algebra $\mI_n$, the factor algebra $\mI_n/\gp$ is isomorphic to the algebra $\mI_{n,m}$ for some $m= m(\gp )$. Since the element $\der_i$ is a left regular element of the ring $\mI_1(i)$ and $B_1(i) = \mI_1(i) / F(i)$. The elements $\der_1, \ldots , \der_n$ are left regular elements of the ring $\mI_{n,m}$ for all $m=0,1, \ldots , n$. Now, by (\ref{IIq1}), the elements $\der_1, \ldots , \der_n$ are left regular elements of the algebra $A$.

3. Statement 2 follows from statement 1 in view of the involution $*$ and the fact that $\der_1^*=\int_1, \ldots , \der_n^*= \int_n$.

4. (i) $S${\em is a left Ore set of} $A$: Recall that $S_{\der_1, \ldots , \der_n}^{-1}\mI_n\simeq B_n$ and $\ass_{\mI_n}(S_{\der_1, \ldots , \der_n})=\ga_n$. Clearly, $I\subseteq \ga_n$ since the ideal $\ga_n$ is the largest proper ideal of the algebra $\mI_n$. The multiplicative set $S$ is a left Ore set of $\mI_n$. Hence, its image in the factor algebra $A= \mI_n/I$ is a left Ore set provided $S\cap I=\emptyset$ but this is obvious ($\emptyset = S_{\der_1, \ldots , \der_n}\cap \ga_n \supseteq S\cap I$).

(ii) $S$ {\em is a left denominator set of} $A$: The statement (ii) follows from the statement (i) and statement 2.

(iii) $T$ {\em is a right denominator set of} $A$: The statement (iii) follows from the statement (ii), statement 1  and the fact that  $T= S^*$.

Let $\gp = \gp_{i_1}+\cdots + \gp_{i_m}$.

(iv) $\ass_{l, A}(S)=(\gp +I)/I$ {\em and} $S^{-1}A\simeq A/ \ass_{l, A}(S)\simeq \mI_n/(\gp +I)$: The short exact sequence of $\mI_n$-modules $0\ra I\ra \mI_n \ra A\ra 0$ yields the short exact sequence of $S^{-1}\mI_n$-modules $0\ra S^{-1}I\ra S^{-1}\mI_n \stackrel{\v}{\ra } S^{-1}A\ra 0$. By the statement (ii), $S^{-1}A$ is an algebra and $\v$ is an algebra homomorphism (since $\v (s^{-1}a)= s^{-1}a$). Since $S^{-1}\mI_n= \mI_n/ \gp$ (the operation of localization at $S$ on the left  is equal to the operation of taking modulo the ideal $\gp$) , $S^{-1} A\simeq \mI_n / (\gp + I)$ and $\ass_{l, A}(S)=(\gp +I)/I$.

(v) $\ass_{r, A}(T)=(\gp +I)/I$ {\em and} $AT^{-1}\simeq A/ \ass_{r, A}(T)\simeq \mI_n/(\gp +I)\simeq S^{-1}A$: Since $T^*=S$, we see that $AT^{-1}=(S^{-1}A)^*\simeq (\mI_n / (\gp +I))^*= \mI_n / (\gp^*+I^*)=  \mI_n / (\gp+I)\simeq =S^{-1}A$ since all the ideals of the algebra $\mI_n$ are $*$-invariant; and
$$ \ass_{r, A}(T) =(\ass_{r, A}(T))^{**} =(\ass_{l, A^*}(T^*))^* =(\ass_{l, A}(S) )^* \stackrel{{\rm st.}\, 1}{=} \ass_{l, A}(S) \stackrel{(iv)}{=}(\gp +I)/I.\;\;\; \Box$$

A ring homomorphism $A\ra B$ is called an {\em extension} of $A$. An extension $A\ra B$ is called a {\em left} (resp., {\em right}) {\em flat} if ${}_AB$ (resp., $B_A$) is a flat $A$-module. An extension is called a {\em flat extension} if it is left and right flat. If $S$ is a left denominator set of $A$ then the extension $A\ra S^{-1}A$ is right flat. If $T$ is a right  denominator set of $A$ then the extension $A\ra T^{-1}A$ is left flat. If, in addition, $S^{-1}A\simeq AT^{-1}$, then the extension $A\ra S^{-1}A$ is flat. In particular, if $S$ is a left and right denominator set of $A$ then the extension $A\ra S^{-1}A\simeq AS^{-1}$ is flat.  A left (resp., right; left and right) extension $A\ra B$ is called a {\em left} (resp., {\em right; left and right}) {\em faithfully flat extension} if the functor $-\t_AB$ (resp., $ B\t_A-$; $-\t_AB$ and $B\t_A-$) is faithful, i.e., a nonzero module is mapped in a nonzero module.

\begin{proposition}\label{A1May17}
Let $I$ be a nonzero  ideal of the algebra $\mI_n$, $A=\mI_n/I$ and $A_i= S_{\der_i}^{-1}A\simeq AS_{\int_i}^{-1}$ for $i=1, \ldots , n$ (Proposition \ref{A30Apr17}.(4)). Then
\begin{enumerate}
\item The set $\{ A\ra A_i\}_{i=1, \ldots , n}$ is a set of (left and right) flat extensions such that the extension $A\ra A^{ff}:=\prod_{i=1}^n A_i$ is a  (left and right) faithfully  flat extension.
\item For each $i=1, \ldots , n$, $\mI_n = \mI_1(i) \t \mI_{n-1}[i]$, where $\mI_{n-1}[i]:= \bigotimes_{j\neq i}\mI_1(i)$, and $A_i\simeq B_1(i)\t (\mI_{n-1}[i]/I_i)$ for some ideal $I_i$ of $\mI_{n-1}[i]$.
\end{enumerate}
\end{proposition}

{\it Proof}. 1. Since $A_i=S_{\der_i}^{-1}A\simeq AS_{\int_i}^{-1}$ (Proposition \ref{A30Apr17}.(4)),  the set $\{ A\ra A_i\}_{i=1,\ldots , n}$ is a set of (left and right) flat extensions.

Since $A^*=A$ and $A_i^*\simeq A_i$ for $i=1, \ldots , n$,  it suffices to show that the extension
$A\ra A^{ff}$ is right faithfully flat, i.e., if $A^{ff}\t_AM=0$ for some $A$-module $M$ then $M=0$. Clearly, $A^{ff}\t_AM=0$ iff $S_{\der_i}^{-1} M=0$ for all $i=1, \ldots , m$. Suppose that $M\neq 0$, then we can choose a nonzero element $m$ of $M$ such that $\der_1m =0 , \ldots , \der_m m=0$.  Since the simple $\mI_n$-module $P_n$ is isomorphic to the $\mI_n$-module $\mI_n / \mI_n (\der_1, \ldots , \der_n)$ (\cite[Proposition 3.8]{algintdif}), the $\mI_n$-module $\mI_n m$ is isomorphic to the $\mI_n$-module $P_n$. Since $M$ is an $A$-module, it is also an $\mI_n$-module such that $IM=0$. In particular, $I\mI_n m =0$ but the $\mI_n$-module $P_n$ is faithful (\cite[Proposition 3.8]{algintdif}), a contradiction. Therefore, $M=0$.

2. The algebra $B_1(i)$ is a central simple algebra and the algebra $A_i$ is an epimorphic image of the algebra $B_1(i)\t \mI_{n-1}[i]$. Therefore, $A_i\simeq B_1(i)\t (\mI_{n-1}[i]/I_i)$ for some ideal $I_i$ of the algebra $\mI_{n-1}[i]$.  $\Box $

Let $I$ be an nonzero  ideal of the algebra $\mI_n$, By Theorem \ref{b10Oct9}.(8), each minimal prime $\gq$ over $I$ is a unique sum $\gq = \sum_{i\in I(\gq )}\gp_i$ for a unique non-empty subset $I(\gq )$ of $\{ 1, \ldots , n\}$ and the set $\{ I(\gq ) \, | \, \gq \in \min (I)\}$ consists of  incomparable elements (for all distinct $\gq , \gq'\in \min (I)$, $\gq \not\subseteq \gq'$ and $\gq \not\supseteq \gq'$). Let $S(\gq )$ (resp., $T(\gq )$) be the monoid generated by the elements $\{ \der_i \, | \, i\in I(\gq ) \}$ (resp., $\{ \int_i \, | \, i\in I(\gq ) \}$). By Proposition \ref{A30Apr17}.(4), the set $S(\gq )$ (resp., $T(\gq )$) is a left (resp., right) denominator set of the algebras $\mI_n$ and $A$ such that $S(\gq )^* = T(\gq )$, $\ass_{l, A}S(\gq )) = \gq/I = \ass_{r, A}(T(\gq ))$ and $S(\gq )^{-1}A\simeq AT(\gq )^{-1} \simeq A/\gq $. Therefore, the extension (\ref{IIq1}) is a (left and right) flat.

\begin{proposition}\label{B1May17}
Let $I$ be a nonzero  ideal of the algebra $\mI_n$ and $A=\mI_n/I$. Then the extension
$A\ra \widetilde{A}=\prod_{\gq \in \min (I)} \mI_n/  \gq$ is a (left and right) faithfully flat extension.
\end{proposition}

{\it Proof}. We have seen above that the extension is (left and right) flat. Since $f(a^*) = f(a)^*$ for all elements $a\in A$, to finish the proof it suffices to show that the extension is left faithful. By \cite[Proposirion 2.3]{MR}, we have to show that $J\widetilde{A}\neq \widetilde{A}$ for all proper right ideals $J$ of $A$. Suppose that $J\widetilde{A}=\widetilde{A}$ for some proper right ideal $J$ of $A$, we seek a contradiction. Then for all ideals $\gq \in \min (I)$,  $J(\mI_n/\gq ) = \mI_n/\gq$, and so $J+\gq = \mI_n$ for all ideals $\gq \in \min (I)$. For each ideal $\gq \in \min (I)$, there are elements $j(\gq) \in J$ and $ q(\gq ) \in \gq$ such that $j(\gq ) + q(\gq ) =1$. In $A$, $0=\prod_{\gq \in \min (I)} q(\gq ) = \prod_{\gq \in \min (I)}(1-j(\gq )) = 1-j$ for some element $j\in J$ (the order of the elements in the products above is an arbitrary but fixed) and so $1=j\in J$, a contradiction. $\Box$

\begin{theorem}\label{T1.2THM}
{\rm (\cite[Theorem 1.2]{THM}.)} Suppose that $\{A\ra A_{\alpha  }\,|\,\alpha  \in I \}$ is a set of flat extensions of a ring $A$ such that  $\wdim (A)<\infty ,$ $\nu =\sup \{ \wdim(A_{\alpha  })\,|\,\alpha  \in I\}<\infty ,$ $\mu =\sup \{\fd ({}_AM)\,|\,M $ is a cyclic $A$-module such that $A_{\alpha  }\t_AM=0$ for all $\alpha  \in I\};$ the set $I$ is either a finite set  or the ring $A$ is right coherent. Then
\begin{enumerate}
\item either  $\wdim (A)=\mu $ or $\mu < \wdim (A) \leq \nu.$

\item If, additionally,  $\nu \leq  \wdim (A)$ (for example, all the rings $A_{\alpha  }$ are two-sided localizations of $A$), then  $\wdim (A)=\max \{\mu ,\nu \}$.
\end{enumerate}
\end{theorem}

Theorem \ref{1May17} and Theorem \ref{BGanya1May17} are generalizations of Theorem \ref{T1.2THM}. They are used in the proof of Theorem \ref{29Apr17}.

\begin{theorem}\label{1May17}
Let $A$ be a ring, $S_1, \ldots S_n$ be left denominator sets of the ring $A$ such that the rings $A_i:= S_i^{-1}A$ ($i=1, \ldots , n$) are flat left $A$-modules. Then $\wdim (A) = \max \{ \nu , \mu \}$ where $\nu = \max \{ \wdim (A_1) , \ldots , \wdim (A_n)\}$ and $\mu =\sup \{ \fd({}_AL)\, | \, L$ is a cyclic $S_i$-torsion $A$-module for $i=1, \ldots , n\}$.
\end{theorem}

{\it Proof}. To prove the theorem we use induction on $n$. The case $n=1$ follows from \cite[Theorem 2.5]{glgwa}. So, let $n>1$ and we assume that the equality holds for all $n'<n$. Clearly, $\nu \leq \max \{ \nu , \mu \} \leq \wdim (A)$. So, if $\max \{ \nu , \mu \}=\infty$ there is nothing to prove. So, we assume that $\max \{ \nu , \mu \}<\infty$. It suffices to show that $\wdim (A)<\infty$ since then the equality would follows from  Theorem \ref{T1.2THM}.(2).
By induction on $n$, $\wdim (A) = \sup \{ \nu_{n-1}, \mu_{n-1}\}$ where $\nu_{n-1}= \max \{ \wdim (A_1) , \ldots , \wdim (A_{n-1})\}$ and $\mu_{n-1} =\sup \{ \fd({}_AL)\, | \, L$ is a cyclic $S_i$-torsion $A$-module for $i=1, \ldots , n-1\}$. Let $L$ be an $S_i$-torsion left $A$-module for $i=1, \ldots , n-1$, and $L_0=\tor_{S_n}(L)$. We have short exact sequences of $A$-modules
\begin{equation}\label{LLL1}
0\ra L_0\ra L\ra \bL \ra 0,
\end{equation}
\begin{equation}\label{LLL2}
0\ra \bL \ra S_n^{-1}L\ra L_1 \ra 0.
\end{equation}
Since $A\ra A_n$ is a flat extension, $\fd({}_AS_n^{-1}L)\leq \fd ({}_{A_n}S_n^{-1}L)\leq \wdim (A_n)<\infty$ (since for all modules $M_A$ and ${}_{A_n}N$ and $m\geq 0$, ${\rm Tor}^A_m(M,N)={\rm Tor}^A_m(M,A_n\t_{A_n}N)={\rm Tor}^{A_n}_m(M\t_AA_n,N)$). Clearly, the $A$-modules $L_0$ and $L_1$ are $S_i$-torsion for $i=1, \ldots , n$. Therefore, $\fd (L_j)\leq \mu < \infty$ for $j=0,1$. By (\ref{LLL2}),
$$\fd ({}_A\bL ) \leq \max \{ \fd({}_AS_n^{-1}L), \fd({}_AL_1)\}\leq \max \{ \wdim (A_n) , \mu\} <\infty .$$ Then, by (\ref{LLL1}), $\fd ({}_AL)\leq \max \{ \fd ({}_AL_0), \fd ({}_A\bL )\} \leq \max \{ \wdim (A_n), \mu \} <\infty$, as required.  $\Box $

\begin{theorem}\label{BGanya1May17}
Let $A$ be a ring. For each  $i=1, \ldots , n$, $S_i$ be left denominator set and $T_i$ be a right denominator set of $A$ such that $A_i:= S_i^{-1}A\simeq AT_i^{-1}$.  Then $\wdim (A) = \max \{ \nu , \mu \} = \max \{ \nu , \mu' \}$ where $\nu = \max \{ \wdim (A_1) , \ldots , \wdim (A_n)\}$, $\mu =\sup \{ \fd({}_AL)\, | \, L$ is a cyclic $S_i$-torsion left  $A$-module for $i=1, \ldots , n\}$ and $\mu' =\sup \{ \fd(N_A)\, | \, N$ is a cyclic $T_i$-torsion right  $A$-module for $i=1, \ldots , n\}$. If, in addition, the flat extension $A\ra \prod_{i=1}^nA_i$ is a left or right faithfully flat extension then $\wdim (A) = \nu$.
\end{theorem}

{\it Proof}. The theorem is a particular case of Theorem \ref{1May17}, hence $\wdim (A) = \max \{ \nu , \mu \}$. Then the second equality, $\wdim (A) = \max \{ \nu , \mu' \}$, follows from the first one in view of symmetry of the weak global dimension with respect left and right modules. If, in addition, the flat extension $A\ra \prod_{i=1}^nA_i$ is a left (resp., right) faithfully flat extension then $\mu' =0$ (resp.,  $\mu =0$) and so  $\wdim (A) = \nu$.  $\Box$

\begin{theorem}\label{29Apr17}
\begin{enumerate}
\item Let $\CA$ be a factor algebra of the algebra $\mI_{n,m}=B_{n-m}\t \mI_m$ where $0\leq m \leq n$. Then $\wdim (\CA) = n$.
\item Let $A$ be a factor algebra of $\mI_n$. Then $\wdim (A) = n$.
\end{enumerate}
\end{theorem}

{\it Proof}. 1. To prove the equality we use induction on $n$. For $n=1$, there are only two options for the algebra $\CA$. Namely, either $A= \mI_1$ or $A= B_1$. By Theorem \cite[Theorem 6.2]{algintdif}, $\wdim (\mI_1) =1$. It is a well-known fact that the algebra $B_1$ is a Noetherian algebra of global dimension 1. In particular, $\wdim (B_1)=1$. Suppose that $n>1$ and the result holds for all $n'<n$. Now, we use the second induction on $m=0,1,\ldots , n$. The initial case $m=0$ is obvious since the algebra $\mI_{n,0}=B_n$ is a simple algebra of global dimension $n$. So, we assume that $m>0$ and the result holds for all $m'$ such that $0\leq m'<m$. The algebras $B_k$ ($k\geq 0$) are central simple algebras. So, the factor algebra $\CA$ of the algebra $\mI_{n,m}= B_{n-m, m}\t \mI_m$ is isomorphic to the algebra $B_{n-m}\t (\mI_m / I_m)$ for some ideal $I_m$ of $\mI_m$. If $I_m=0$ then $\CA = \mI_{n,m}$ and $\wdim (\CA) = \wdim (\mI_{n,m})=n$ (\cite[Theorem  6.2]{algintdif}). We can assume that $I_m\neq 0$. Then $I_m$ contains the least nonzero ideal $F^{\t m}$ of the algebra $\mI_m$ (Theorem \ref{b10Oct9}.(8)). By Proposition \ref{A30Apr17}.(4), for each $i=n-m+1, \ldots , n$, consider  the ring extension
\begin{equation}\label{CAAi}
\CA \ra \CA_i:= S_{\der_i}^{-1}\CA \simeq \CA S_{\int_i}^{-1}\simeq B_{n-m}\t B_1\t (\mI_{m-1}/J_i)
\end{equation}
 for some ideal $J_i$ of the algebra $\mI_{m-1}$. By induction, $\wdim (\CA_i) = n$ for all $i=n-m+1, \ldots , n$. By Theorem \ref{BGanya1May17}, $\wdim (\CA ) = \max \{ \nu , \mu \}$ where $\nu = \max \{ \wdim (\CA_i) \, | \, i=n-m+1, \ldots , n\}=n$ and $\mu = \sup \{ \fd ({}_\CA L) \, | \, L\in \mE\}=0$ since $\mE :=\{ L \, | \, L$ is a cyclic $S_i$-torsion left  $\CA$-module for $i=n-m+1, \ldots , n\}=\{ 0\}$, by Proposition \ref{A1May17} as $I_m\neq 0$. Therefore, $\wdim (\CA ) =n $.

2. Statement 2 is a particular case of statement 1 since $\mI_n = \mI_{n,n}$.  $\Box $

Theorem \ref{3May17} is an analogue of Hilbert's Syzygy Theorem for all factor algebras of the algebra $\mI_n$ in case of the weak global dimension.

\begin{theorem}\label{3May17}
Let $K$ be an algebraically closed uncountable field of
characteristic zero, $\CA$ be a factor algebra of $\mI_{n,m}=B_{n-m}\t \mI_m$, where $0\leq m \leq n$,  and $B$ be a left Noetherian finitely
generated algebra over $K$. Then
\begin{enumerate}
\item  $\wdim (\CA\t B) = \wdim (\CA) +\wdim (B)  = n+\wdim
(B)$.
\item  $\wdim (A\t B) = \wdim (A) +\wdim (B)  = n+\wdim (B)$ for all factor algebras $A$ of $\mI_n$.
\end{enumerate}
\end{theorem}

{\it Proof}. 1.  To prove the theorem we use induction on $n$. For $n=1$, there are only two options for the algebra $\CA$. Namely, either $\CA= \mI_1$ or $\CA = B_1$. In both cases, the result is a particular case of \cite[Theorem 6.5]{algintdif}. Suppose that  $n>1$ and we assume that the equality holds for all $n'<n$. Now, we use the second induction on $m=0,1,\dots , n$. In the initial case when $m=0$, $\CA = B_n$ (since $\mI_{n,0}=B_n$ is a simple Noetherian finitely generated algebra) and the result is known (see \cite[Theorem 6.5]{algintdif}). So, we may assume that $m>0$ and the equality holds for all $m'<m$. Notice that $\mI_{n,m}=B_{n-m}\t \mI_m$ and the algebra $B_{n-m}$ is a central simple algebra. Therefore, $\CA \simeq B_{n-m}\t (\mI_m/I_m)$ for some ideal $I_m$ of the algebra $\mI_m$. If $I_m=0$ then the result is precisely \cite[Theorem 6.5]{algintdif}.  So, we assume that $I_m\neq 0$. By Proposition \ref{A30Apr17}.(4), for each $i=n-m+1, \ldots  , n$, consider the extension $\CA \ra \CA_i$, see (\ref{CAAi}). By Proposition \ref{A1May17}.(1), the extension $\CA \ra \widetilde{A}:=\prod_{i=n-m+1}^n \CA_i$ is faithfully flat, hence so is the extension $A\t B\ra \widetilde{A}\t B$. By (\ref{CAAi}) and the induction on $m$, for all $i$,
 $$\wdim (\CA_i\t B) = \wdim (\CA_i) +\wdim (B)  = n+\wdim
(B)$$ (since $\wdim (\CA_i) =n$, by Theorem \ref{29Apr17}). Now, by Theorem \ref{BGanya1May17} and Theorem \ref{29Apr17},
 $\wdim (\CA \t B) = \max\{ \wdim (\CA_1\t B)\, | \, i=n-m+1, \ldots , n\}= n+\wdim (B)= \wdim (\CA) + \wdim (B).$

 2. Statement 2 is a particular case of statement 1 since $\mI_n = \mI_{n,n}$. $\Box$


\section{The weak  global dimension of factor algebras of
$\mA_n$}\label{WDIMFAN}

The aim of this section is to show that the weak global dimension of all factor algebras of the Jacobian algebra $\mA_n$ is $n$ (Theorem \ref{An29Apr17}) and that an analogue of Hilbert's Syzygy Theorem holds for them (Theorem \ref{A3May17}).

A classification of ideals of the Jacobian algebra $\mA_n$ is given in \cite{Bav-Jacalg}.

{\it Definition}. Let $A$ and $B$ be algebras, and let $\CJ (A)$
and $\CJ (B)$ be their lattices of ideals. We say that a bijection
$f: \CJ (A) \ra \CJ (B)$ is an {\em isomorphism} if $f(\ga *\gb )
= f(\ga )*f(\gb )$  for $*\in \{ +, \cdot, \cap \}$, and in this
case we say that the algebras $A$ and $B$ are {\em ideal
equivalent}. The ideal equivalence is an equivalence relation on
the class of algebras.

The next theorem shows that the algebras $\mA_n$ and $\mI_n$ are
ideal equivalent.
\begin{theorem}\label{7Oct9}
{\rm (\cite[Theorem 3.1]{Bav-Jacalg}.)} The restriction map $\CJ ( \mA_n) \ra \CJ (\mI_n)$, $\ga \mapsto
\ga^r:= \ga \cap \mI_n$, is an isomorphism (i.e.,  $(\ga_1*\ga_2)^r=
\ga_1^r*\ga_2^r$ for $*\in \{ +, \cdot, \cap \}$) and its inverse
is the extension map $\gb\mapsto \gb^e:= \mA_n \gb \mA_n$.
\end{theorem}
 It follows from the explicit description of ideals of the algebra $\mI_n$ (Theorem \ref{b10Oct9}.(9)) and $\mA_n$ (\cite[Theorem 3.1]{Bav-Jacalg}) that the lattice isomorphism
\begin{equation}\label{JIJA}
\CJ ( \mI_n) \ra \CJ (\mA_n), \;\; \gb \mapsto
\gb^e:= \mA_n\gb \mA_n= S^{-1}\gb = \gb S^{-1},
\end{equation}
can be written via the localizations at $S$, see (\ref{mAnlIn}).

\begin{theorem}\label{An29Apr17}
\begin{enumerate}
\item Let $\CA$ be a factor algebra of the algebra $\mA_{n,m}=\CA_{n-m}\t \mA_m$ where $0\leq m \leq n$. Then $\wdim (\CA) = n$.
\item Let $A$ be a factor algebra of $\mA_n$. Then $\wdim (A) = n$.
\end{enumerate}
\end{theorem}

{\it Proof}. 1. Since the algebra $\CA_{n-m}$ is a central simple algebra, $\CA \simeq \CA_{n-m}\t (\mA_m/ J)$ for some ideal $J$ of the algebra $\mA_m$. By (\ref{mAny}),
\begin{equation}\label{SxCA}
S^{-1}_{y_{n-m+1}, \ldots y_n}\CA\simeq \CA S^{-1}_{x_{n-m+1}, \ldots x_n}\simeq \CA_{n-m}\t \CA_m\simeq \CA_n
\end{equation}
since $\ass_{l,\mA_m}(S_{y_{n-m+1}, \ldots y_n})=\ass_{r,\mA_m}(S_{x_{n-m+1}, \ldots x_n})=\ga_m$ is the largest proper ideal of the algebra $\mA_m$ (it follows from the inclusion $J\subseteq \ga_m$ that $S^{-1}_{y_{n-m+1}, \ldots y_n} (\mA_m/ J)\simeq S^{-1}_{y_{n-m+1}, \ldots y_n} \mA_m \simeq \CA_m$ and $(\mA_m/ J)S^{-1}_{x_{n-m+1}, \ldots x_n}\simeq \mA_mS^{-1}_{x_{n-m+1}, \ldots x_n}\simeq \CA_m$). By Theorem \ref{bigAn}.(1,3) and (\ref{SxCA}),
$$ n = \gldim (\CA_n) = \wdim (\CA_n)= \wdim (S^{-1}_{y_{n-m+1}, \ldots y_n} \CA) \leq \wdim (\CA).$$
So, by (\ref{JIJA}), $\mA_m/J\simeq S^{-1}(\mI_m/I)$ for some
ideal $I$ of $\mI_m$. So,
\begin{equation}\label{ASA}
\CA \simeq \CA_{n-m}\t (\mA_m / J)\simeq (S^{-1}_{y_1, \ldots , y_{n-m}}\mI_{n-m})\t S^{-1}(\mI_m/ I).
\end{equation}
Now,
$$\wdim (\CA) \stackrel{(\ref{ASA})}{\leq} \wdim (\mI_{n-m}\t (\mI_m/I))=n, \;\; {\rm (by \; Theorem \; \ref{29Apr17}.(1)).}$$
Therefore, $\wdim (\CA ) =n$.

2.  Statement 2 is a particular case of statement 1 when $m=n$ since $\mA_{n,n}= \mA_n$.  $\Box$

\begin{theorem}\label{A3May17}
Let $K$ be an algebraically closed uncountable field of
characteristic zero, $\CA$ be a factor algebra of $\mA_{n,m}=\CA_{n-m}\t \mA_m$, where $0\leq m \leq n$,  and $B$ be a left Noetherian finitely
generated algebra over $K$. Then
\begin{enumerate}
\item  $\wdim (\CA\t B) = \wdim (\CA) +\wdim (B)  = n+\wdim
(B)$.
\item  $\wdim (A\t B) = \wdim (A) +\wdim (B)  = n+\wdim (B)$ for all factor algebras $A$ of $\mA_n$.
\end{enumerate}
\end{theorem}

{\it Proof}. 1.  The algebra $B$ is a left Noetherian algebra, hence so is the algebra $\CA_n \t B$. By Theorem \ref{bigAn}.(5),
$$ \wdim (\CA_n \t B) =\lgldim (\CA_n \t B) = \lgldim (\CA_n) +\lgldim (B) = n+\wdim (B). $$
We keep the notation of the proof of Theorem \ref{An29Apr17}.
\begin{eqnarray*}
\wdim (\CA_n\t B)  &\stackrel{(\ref{SxCA})}{=}& \wdim (S^{-1}_{y_{n-m+1}, \ldots , y_n}\CA \t B) \leq \wdim (\CA \t B)\stackrel{(\ref{ASA})}{\leq}\wdim(\mI_{n-m}\t\mI_m/I\t B) \\
 &=&  n+\wdim (B), \;\; {\rm (by\;  Theorem \; \ref{3May17}.(2))}\\
 &=& \wdim (\CA ) + \wdim (B),  \;\; {\rm (by\;  Theorem \; \ref{An29Apr17}.(2))}.
\end{eqnarray*}

2. Statement 2 is a particular case of statement 1 when $m=n$ since $\mA_{n,n}=\mA_n$. $\Box$

Department of Pure Mathematics

University of Sheffield

Hicks Building

Sheffield S3 7RH

UK

email: v.bavula@sheffield.ac.uk

\end{document}